\documentclass{article}
\usepackage{graphicx} 

\usepackage[english]{babel}

\usepackage[letterpaper,top=2cm,bottom=2cm,left=3cm,right=3cm,marginparwidth=1.75cm]{geometry}

\usepackage[utf8]{inputenc}
\usepackage{authblk,amsmath,amsfonts,amssymb,amsthm}
\usepackage{graphicx}
\usepackage{longtable}
\usepackage{epsfig}
\usepackage{tikz}
\usepackage{tikz-cd}
\usepackage{nicefrac}
\usepackage{csquotes}
\usepackage{mathrsfs}
\usepackage{pdfpages}
\usepackage{framed}
\usepackage{placeins}
\usepackage{fullpage}
\usepackage{graphicx}
\usepackage{mathtools}
\usepackage{wrapfig}
\usepackage{booktabs}
\usepackage{stmaryrd}
\usepackage{color}
\usepackage{float}
\usepackage{fancyhdr}
\usepackage{eurosym}
\usetikzlibrary{patterns}
\usepackage{todonotes}
\usepackage{hyperref}
\usepackage{nameref}
\newcommand{\refwithname}[1]{\ref{#1} (\nameref{#1})}
\newcommand{\refastitle}[1]{\nameref{#1} {\normalfont(c.f. lemma \ref{#1})}}
\usepackage{subcaption}
\usepackage[
  deletedmarkup=sout,
  authormarkup=none,
  todonotes]{changes}
  \newcommand{\stkout}[1]{\ifmmode\text{\sout{\ensuremath{#1}}}\else\sout{#1}\fi}

\newcommand{\jump}[1]{\left\llbracket #1 \right\rrbracket}

\newcommand{\cT}{\mathcal{T}}

\newcommand{\cI}{\mathcal{I}}
\newcommand{\cF}{\mathcal{F}}
\newcommand{\Ctr}[1][X]{C_{tr%
  \IfStrEq{#1}{X}{}{,#1}%
}}

\newcounter{genC}
\newcommand{\genC}[1][\thegenC]{C}

\newcommand{\cFwoo}{\mathcal{F}_{\setminus o}}   
\newcommand{\cFwoio}{\mathcal{F}_{\setminus io}}  

\newcommand{\norm}[1]{\left\lVert #1 \right\rVert}
\newcommand{\abs}[1]{\left\lvert #1 \right\rvert}

\newcommand{\seminormbeta}[1]{\left\lvert #1 \right\rvert_{\beta}}
\newcommand{\normuwb}[1]{\left\lVert #1 \right\rVert_{\text{uwb}}}
\newcommand{\normuwbStern}[1]{\left\lVert #1 \right\rVert_{\text{uwb},*}}
\newcommand{\VhStern}{\mathcal{V}_{*}^0}

\newcommand{\ahDoD}{a_h^{{DoD}}}
\newcommand{\AhDoD}{A_h^{{DoD}}}
\newcommand{\Ahupw}{A_h^{{upw}}}
\newcommand{\average}[1]{\left\lbrace\mskip-6mu\left\lbrace #1 \right\rbrace\mskip-6mu\right\rbrace}

\newcommand{\dt}{\Delta t}
\newcommand{\xipin}{\xi_{\pi}^n}

\newcommand{\xihn}{\xi_h^n}
\newcommand{\xihnplus}{\xi_h^{n+1}}
\newcommand{\trace}[1]{\operatorname{tr}(#1)}
\newcommand{\betamean}[2]{\overline{#1\rule{0pt}{1.8ex}}^{{\beta},#2}}
\newcommand{\vweighted}{\betamean{v|_E}{e}}
\newcommand{\taueps}{\tau_{\varepsilon}}
\newcommand{\epsinv}{{\textstyle \frac 1 \varepsilon}}

\newcommand{\Sandra}[1]{\textcolor{green!50!blue!80!white}{#1}}

\newcommand{\Gunnar}[1]{\textcolor{green!10!blue!70!red!100}{#1}}

\usepackage{todonotes}
\usepackage{pdfpages}
\usepackage{enumitem}
\newlist{todolist}{itemize}{2}
\setlist[todolist]{label=$\square$}
\usepackage{pifont}

\colorlet{MyColorBlue}{blue!50}

\newtheorem{theorem}{Theorem}[section]

\newtheorem{lemma}[theorem]{Lemma}

\theoremstyle{definition}
\newtheorem{definition}[theorem]{Definition}
\newtheorem{notation}[theorem]{Notation}
\theoremstyle{remark}
\newtheorem{remark}[theorem]{Remark}

\title{Error analysis of a first-order DoD cut cell method\\ for 2D unsteady advection}
\author{Gunnar Birke\thanks{\url{gunnar.birke@uni-muenster.de}, \url{christian.engwer@uni-muenster.de}, Applied Mathematics M\"unster: Institute for Analysis and Numerics, University of M\"unster, Einsteinstra{\ss}e 62, 48149 Münster, Germany} , Christian Engwer$^*$, Jan Giesselmann\thanks{\url{jan.giesselmann@tu-darmstadt.de}, Department of Mathematics, Technical University of Darmstadt,  Dolivostr 15, 64293 Darmstadt, Germany
} , Sandra May\thanks{\url{may@math.tu-berlin.de}, Institute of Mathematics, TU Berlin, Stra{\ss}e des 17. Juni 136,
10623 Berlin, Germany}\:\:\thanks{Department of Information Technology, Uppsala University, L\"agerhyddsv\"agen 1, 751 05 Uppsala, Sweden} }

\date{}

\begin{document}

%
%

\maketitle

\begin{abstract}
In this work we present an a priori error analysis for solving the unsteady advection equation on cut cell meshes along a straight ramp in two dimensions. The space discretization uses a lowest order upwind-type discontinuous Galerkin scheme involving a \textit{Domain of Dependence} (DoD) stabilization to correct the update in the neighborhood of small cut cells. Thereby, it is possible to employ explicit time stepping schemes with a time step length that is independent of the size of the very small cut cells.

Our error analysis is based on a general framework for error estimates for first-order linear partial differential equations that relies on consistency, boundedness, and discrete dissipation of  the discrete bilinear form. We prove these properties for the space discretization involving DoD stabilization. This allows us to prove, for the fully discrete scheme, a quasi-optimal error estimate of order one half in a norm that combines the $L^\infty$-in-time $L^2$-in-space norm and a seminorm that contains velocity weighted jumps. We also provide corresponding numerical results.
\end{abstract}

\noindent
\textbf{Keywords:} cut cell, discontinuous Galerkin method, DoD Stabilization, a priori error estimate, unsteady advection

\noindent
\textbf{AMS subject classifications:} 65M60, 65M15, 65M20, 35L02\\

\section{Introduction}

So called \textit{cut cell} or \textit{embedded boundary} meshes have become popular in recent years due to the easyness of the mesh generation process: a given geometry is simply cut out of a structured, often Cartesian, background mesh. This results in \textit{cut cells} where the object intersects the background mesh. Cut cells can have various shapes, the sizes of neighboring elements can differ by several orders of magnitudes, and (most difficult to deal with) can become arbitrarily small. For the solution of time-dependent hyperbolic conservation laws this causes the \textit{small cell problem}: standard explicit time stepping schemes are not stable on small cut cells when the time step length is chosen with respect to the background mesh.

One way to address this issue is by using \textit{cell merging/cell agglomeration}, see, e.g., \cite{Bayyuk_Powell_vanLeer,Quirk1994,QIN201324,Oberlack2016}. In this approach small cut cells are merged with larger neighbors until all small cut cells are gone. The resulting mesh still contains unstructured (cut) cells along the embedded boundary but the issue of smallness is gone.

Another approach is to keep the small cut cells and to stabilize them. Over the years, different methods have been developed in the context of finite volume and discontinuous Galerkin (DG) schemes to solve the small cell problem as described above. Two well-known approaches are the \textit{flux redistribution} method, see, e.g., \cite{Chern_Colella,Colella2006}, and the $h$\textit{-box} method \cite{Berger_Helzel_Leveque_2003,Berger_Helzel_2012}. Newer developments are the \textit{dimensionally split} approach \cite{Klein_cutcell, Gokhale_Nikiforakis_Klein_2018},  
the \textit{mixed explicit implicit} scheme \cite{May_Berger_explimpl}, the extension of the active flux method to cut cell meshes \cite{FVCA_Helzel_Kerkmann}, the \textit{state redistribution} method \cite{Berger_Giuliani_2021,Giuliani_DG}, the \textit{Domain of Dependence (DoD)} stabilization \cite{DoD_2d_linadv_2020,DoD_1d_nonlin_2022}, and the extension of the \textit{ghost penalty} stabilization \cite{Burman2010} to first-order hyperbolic problems on cut cell meshes \cite{Fu_Kreiss_2021, Fu_Frachon_Kreiss_Zahedi_2022}. All of these methods are based on finite volume schemes or stabilize a standard DG scheme. Preliminary work \cite{KaurHicken_DGD_CutCell} in the context of so called \textit{DG difference} schemes indicates that due to the extended support of the underlying basis functions, one might get away without stabilizing the scheme on small cut cells but further investigation is needed.

The derivation of error estimates for these schemes is very challenging. Since stabilized schemes on cut cell meshes keep the resulting potentially
arbitrarily small and very skewed cut cells, any derivation of an error estimate has to handle these cells, leading to many complications. In this work we will present an \textit{error estimate} for the time-dependent linear advection equation in \textit{two dimensions} for the \textit{fully discrete scheme} based on using the \textit{DoD stabilization}. Due to complexity, we will only consider advection along a straight ramp, using piecewise constant functions in space combined with explicit Euler in time. We will show that under a CFL constraint that is independent of the size and geometry of the small cut cells, there holds a result of the form
\begin{equation}
  \| u(t^M,\cdot) - u_h^M \|_{L^2(\Omega)}^2 + \sum_{n=0}^{M-1} \frac{\dt}{4} \seminormbeta{u(t^n,\cdot) - u_h^n}^2
  \leq \mathcal{O}( \dt^2 + h).
 \end{equation}
 The seminorm $\seminormbeta{\cdot}$ will be defined below in \eqref{eq: seminorm beta 2d P0} and consists of face terms. As the numerical results in section \ref{sec: numerical results} will show, the result is sub-optimal in the $L^2$ norm (like many DG-error estimates in two dimensions) but optimal in the seminorm. (Note that the factor $\dt$ on the left hand side essentially cancels with the sum over the number of time steps.)
 This result is very close to \cite[Thm.~3.7]{DiPietro_Ern} where an analogous scheme is studied but on a quasi-uniform and shape regular mesh without DoD stabilization.
 
The proof is based on the powerful proof framework that has been developed and is presented in \cite{BURERNFER2010} and \cite[Ch.~3]{DiPietro_Ern} and the references cited therein, 
and has been adjusted appropriately to the situation of the DoD stabilized scheme for cut cells. 
%
These changes affect the fundamental properties of the bilinear form, i.e., consistency, discrete dissipation, and boundedness, that we will show in section \ref{sec: auxiliary results}.
In particular, dissipation and boundedness are proven with respect to DoD-adjusted seminorms and norms that include an appropriate weighting. This makes their proofs somewhat more involved.  In addition,
new inverse estimates and new projection error estimates had to be proven.
The final proof in section \ref{sec: proof of main result} that puts all the pieces together then follows the original proof except for the adjustment of the consistency error. 
To the best of our knowledge, this is the first error estimate for a fully discrete scheme for solving the unsteady linear advection equation in two dimensions on a cut cell mesh with arbitrarily small cut cells using a stabilized scheme.

The situation in one space dimension is very different. 
For the \textit{one-dimensional} linear advection equation on a model problem, where typically one small cell is inserted in an otherwise equidistant mesh, error estimates have been shown for several schemes. Many use the following trick introduced in \cite{Wendroff_White, Wendroff_White_proc}, which is suitable for finite volume and finite difference schemes: for showing order $p$, in a truncation error analysis of the one step error, one typically finds that the error is only $\mathcal{O}(h^{p})$ (instead of $\mathcal{O}(h^{p+1})$) in the neighborhood of the cut cell. One then defines a new grid solution such that (1) the one step error is $\mathcal{O}(h^{p+1})$ for this new solution and such that (2) the new solution is at each grid point and at each time step only $\mathcal{O}(h^p)$ away from the true solution. In a cut cell context, this approach has first been used to show first- and second-order error estimates for the $h$-box method \cite{Berger_Helzel_Leveque_2003}. It has then been copied to show a first-order result for the dimensionally split approach \cite{Gokhale_Nikiforakis_Klein_2018}, a first-order result for the DoD stabilization \cite{DoD_2d_linadv_2020} (for piecewise constant polynomials, the DG based DoD stabilization and the $h$-box method are essentially the same), and a second-order result for the mixed explicit implicit scheme \cite{May_Laakmann_2024}. Finding such a new grid solution is already challenging in one dimension, and would probably be enormous work in two dimensions.

Again for the time-dependent advection equation in one dimension but in a more \textit{DG}-style approach are the results for the stabilization based on a ghost penalty approach \cite{Fu_Kreiss_2021,Fu_Frachon_Kreiss_Zahedi_2022}: for the semi-discrete setting, using polynomials of degree $p$, order $p+\tfrac 1 2$ is shown for the $L^2$-norm of the spatial error as well as order $p$ when an interface is present.

This contribution is structured as follows: First, in section \ref{sec: section2}, we describe the problem setup, introduce notation, and define the fully discrete problem based on using the DoD stabilization. We also introduce the seminorm $\seminormbeta{\cdot}$ and the norms $\normuwb{\cdot}$ and $\normuwbStern{\cdot}$ that we will use throughout this paper in our proofs. Finally, we collect the major results (concerning consistency, discrete dissipation, and boundedness) that will be necessary to prove the desired result and where the major adjustments were made. 
The proof of these results will be given later in section \ref{sec: auxiliary results}. Before that, in section 
\ref{sec:geom}, we will prove cut cell related estimates such as inverse estimates accounting for the potentially arbitrarily small cut cells and a new result for the projection error. 
In section \ref{sec: proof of main result}, we will then put all the pieces together and prove the main error estimate. In the final section \ref{sec: numerical results} we will present numerical results to examine the optimality of our theoretical results. 


\section{Problem setting and main result}\label{sec: section2}

\begin{figure}
\begin{subfigure}[b]{0.20\textwidth}
  \centering
  \begin{tikzpicture}[scale=1.5]
    \node at (.25,1.25) {\small$\Omega$};
    \draw[semithick] (0,0) rectangle (1.5,1.5);
    \fill[MyColorBlue] (0.3,0.0) --
      (1.5,0.7) --
      (1.51,0.7) --
      (1.51,-0.01) --
      (0.3,-0.01);
    \draw[dashed](0.3,0) -- +(0:0.7) arc (-45:45:0.29) -- cycle;
    \node[] at (0.9,0.15) {$\gamma$} ;  
    \draw[semithick]
    (0.0,0.0) --
    (0.3,0.0) --
    (1.5,0.7) --
    (1.5,1.5) --
    (0.0,1.5) --
    (0.0,0.0);
    \node[] at (0.7,0.7) {$\beta$} ;
    \draw[thick,->] (0.5,0.4)--(1.0,0.69);   
  \end{tikzpicture}
  \caption{Setup of ramp test.}
  \label{fig: ramp geom a}
    \end{subfigure}
  %
  %
  \begin{subfigure}[b]{0.5\textwidth}
  \centering
  \begin{tikzpicture}[scale=1.0]
    \node at (.75,1.75) {\small${\widehat{\mathcal{T}}}_h$};
    \draw[semithick,step=0.25] (0,0) grid (1.5,1.5);
    \node at (1.9,0.75) {\Large$\cap$};
\begin{scope}[xshift=2.3cm]
    \node at (.75,1.75) {\small$\overline{\Omega}$};
    \draw[semithick,fill=black!10!white]
    (0.0,0.0) --
    (0.3,0.0) --
    (1.5,0.7) --
    (1.5,1.5) --
    (0.0,1.5) --
    (0.0,0.0);
    \draw[dashed](0.3,0) -- +(0:0.7) arc (-45:45:0.29) -- cycle;
    \node[] at (0.9,0.15) {$\gamma$} ;
\end{scope}
    \node at (4.4,0.75) {\huge$\rightarrow$};
\begin{scope}[xshift=4.9cm]
    \node at (.75,1.75) {\footnotesize$\mathcal{T}_h$};
    \draw[semithick,step=0.25] (0,0) grid (1.5,1.5);
    \fill[MyColorBlue] (0.3,0.0) --
      (1.5,0.7) --
      (1.51,0.7) --
      (1.51,-0.01) --
      (0.3,-0.01);
    \draw[semithick]
    (0.0,0.0) --
    (0.3,0.0) --
    (1.5,0.7) --
    (1.5,1.5) --
    (0.0,1.5) --
    (0.0,0.0);
\end{scope}
\end{tikzpicture}
\caption{Construction of cut cell mesh.}
\label{fig: ramp geom b}
    \end{subfigure}%
\begin{subfigure}[b]{0.28\textwidth}
  \centering
  \begin{tikzpicture}[scale=1.3]
\draw[] (-0.9,0.5) -- (0.7,0.5);
\draw[] (-0.5,0.9) -- (-0.5,-0.5);
\draw[thick] (-1.0,-0.46) -- (1,0.7);
\fill[MyColorBlue] (-1.0,-0.46) --
      (-1,-0.65) --
      (1,-0.65) --
      (1,0.7) --
      (-1,-0.46);
\draw[->, thick] (-1.2,0.6) -- (-0.7,0.89);
\node[anchor=east] at (-0.9,1.0) {${\beta}$};
\node[anchor=east] at (-0.75,0.125) {\small${E_{in}}$};
\node[anchor=south] at (0.125,0.75) {\small ${E_{out}}$};
\node[anchor=south] at (-0.2,0.1) {\small${E_{cut}}$};
\node[anchor=west] at (0,0) {\small $e_{bdy}$}; 
\node[anchor=west] at (-0.9,0.125) {\small $e_{in}$};
\node[anchor=north] at (0.125,0.75) {\small$e_{out}$};
\end{tikzpicture}
  \caption{Zoom on a triangular cut cell.}
  \label{fig: ramp geom c}
    \end{subfigure}%
  \caption{Geometric considerations and mesh construction for ramp test.}
  \label{fig: ramp 2d cut cell}
\end{figure}
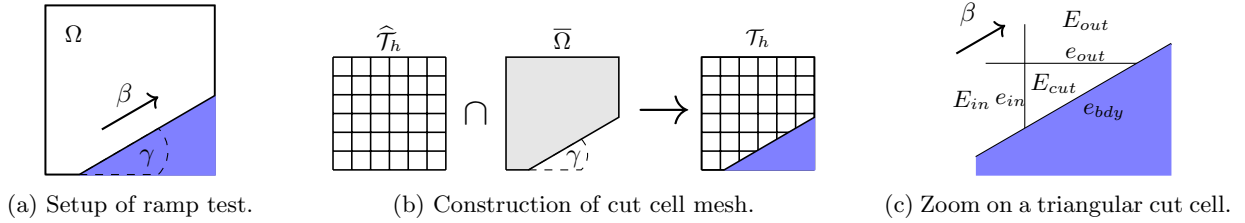

\subsection{Preliminaries}\label{sec: prelim} 
We consider the ramp test setup that is shown in figure \ref{fig: ramp geom a}: the computational domain $\Omega$ consists of a square, where a ramp with angle $\gamma$ has been cut out. We solve the 2d linear advection equation given by
\begin{subequations}\label{eq: lin adv 2d}
\begin{align}
      u_t+\left\langle{\beta},{\nabla u}\right\rangle&=0&\quad&\text{in }(0,T) \times \Omega,\\
  u &=g & &\text{on } (0,T) \times \partial\Omega^{\text{in}},\\
  u &= u_0 & &\text{on }\{ t = 0 \} \times \Omega.
\end{align}
\end{subequations}
Here, $\left\langle \cdot,\cdot \right\rangle$ denotes the Euclidean scalar product in $\mathbb{R}^2$ and $\partial \Omega^{\text{in}} = \{x\in \partial\Omega: \beta(x) \cdot n(x) < 0\}$, where $n \in \mathbb{R}^2$ is the outer unit normal vector on $\partial \Omega$.
The functions for defining inflow boundary conditions and initial data are taken from $g \in H^{1}((0,T) \times \partial\Omega^{\text{in}})$ and $u_0 \in H^1(\Omega)$.
We assume that the velocity field $\beta$ is Lipschitz continuous and satisfies $\nabla \cdot \beta = 0$ and that it is parallel to the ramp, i.e, that along the ramp there holds $\beta \cdot n = 0$. 
For technical reasons, we also assume that the sign of $\beta \cdot n$ is constant along each face of our computational mesh. 
Further, for technical reasons in proving projection error estimates on cut cells in lemma \ref{lem:est_xipi}, we assume that 
$\beta$ does not vanish on the ramp.
{For vector-valued functions, such as $\beta$,} we denote by 
\begin{equation}\label{eq: def beta-inf} 
\| \beta \|_\infty := \max_{x \in \Omega} \abs{\beta(x)}_2, \end{equation}
where $\abs{ \cdot }_2 $ denotes the Euclidean norm on $\mathbb{R}^2$.

For constructing a mesh, we start with a structured background mesh ${\widehat{\mathcal{T}}}_h$, consisting of square elements $\hat{E}$, that discretizes the underlying square domain.
We define the mesh length $h$ as the face length of a Cartesian cell $\hat{E} \in \widehat{\cT}_h$ in $x$- and $y$-direction, i.e., $h_x = h_y = h.$
We intersect ${\widehat{\mathcal{T}}}_h$ with the computational domain $\overline{\Omega}$. This is shown in figure \ref{fig: ramp geom b}. This process creates cut cells $E = \widehat E \cap \overline{\Omega}$ along the ramp. Away from the ramp (and inside $\Omega$), the cells $E$ coincide with the background cells $\widehat{E}$. 
The partition $\cT_h$ is given by
\begin{equation*}
  \cT_h = \left\{~ E = \widehat E \cap \Omega~ \left| ~\widehat E \in \widehat{\mathcal{T}}_h ~\right\}\right..
\end{equation*}
We define the space $\mathcal{V}^0$ of piecewise constant functions on $ \cT_h$ by
\begin{equation}\label{eq: def V_h}
   \mathcal{V}^0 = \left\{ v_h \in L^2(\Omega)  \: \vert \: \forall E \in \cT_h: v_h{\vert_E} \text{ is constant} \right\}.
\end{equation}

\begin{definition}[Space $\VhStern$]
We define 
\[
\VhStern = V + \mathcal{V}^0,
\]
where $V=H^1(\Omega)$. We will assume that the exact solution $u(t,\cdot)$ is an element of $V$. Then $\VhStern$ is the space that contains an error of the form $u(t,\cdot)- w_h$.
\end{definition}

\begin{notation}
In the following we will also use $u(t)$ as a shortcut for $u(t,\cdot)$ {as well as $u^n$ to denote the solution at time $t^n$.}
\end{notation}


Generally, the cut cells that are created along the ramp boundary are either 3-sided, 4-sided, or 5-sided cells. As in \cite{DoD_2d_linadv_2020} we decided to only stabilize 3-sided cut cells. Our analysis will show that this is sufficient to guarantee stability and convergence.
A zoom on such a triangular cut cell $E_{{cut}}$ is shown in figure \ref{fig: ramp geom c}. For the 3 faces there holds that 
\begin{itemize}
\item there is one \textit{inflow} face $e_{{in}}$ with neighboring cell $E_{{in}}$, which is characterized by $\beta \cdot n_{E_{cut}} < 0$,
\item there is one \textit{outflow} face $e_{{out}}$ with neighboring cell $E_{{out}}$, which is characterized by $\beta \cdot n_{E_{cut}} > 0$,
\item and on the \textit{ramp boundary} face $e_{{bdy}}$, there holds $\beta \cdot n_{E_{cut}} = 0$.
\end{itemize}
This notation will be used throughout this work.

\begin{definition}
We denote by $\cI$ {\em the set of stabilized cut-cells}. It is given by
\begin{equation}\label{eq: set of stab cut cells}
\cI = \left\{ E \in \cT_h \: \vert \: E \text{ is triangular cut cell and } \max(\abs{e_{{in}}},\abs{e_{{out}}}) < \tfrac 1 2 h  \right\},
\end{equation}
where $\abs{e}$ denotes the length of a face $e$.
\end{definition}


\begin{notation}[Faces]
We denote by $\cF$ the set of all faces, with $\cF_b$ being the set of boundary faces, $\cF_i$ being the set of interior faces and
$\cF:= \cF_b \cup \cF_i$.
Further, we use the notation
\begin{align*}
\cFwoo &\text{ = the set of all faces that are not outflow faces of stabilized cut cells},\\
\cFwoio & \text{ = the set of faces that are neither inflow nor outflow faces of stabilized cut cells}.
\end{align*}
That is, faces of the kind $e_{{out}}$ are excluded in $\cFwoo$. The set $\cFwoio$ contains neither faces of the kind $e_{{in}}$ nor faces of the kind $e_{{out}}$.
We further define for a cell $E \in \cT_h$ the set of inflow and outflow faces
\begin{align*}
\cF_{in}(E) & = \left\{e\in \cF \: | \: e \in \partial E, \beta \cdot n_E < 0 \right\}  ,\\
\cF_{out}(E) & = \left\{e\in \cF \: | \: e \in \partial E, \beta \cdot n_E > 0 \right\}.
\end{align*}
\end{notation}
For stabilized cut cells $E \in \cI$, $\cF_{in}(E)$ only contains $e_{{in}}$ and $\cF_{out}(E)$ only contains $e_{{out}}$. For any cell $E$ there holds due to the velocity field $\beta$ being incompressible
\begin{equation}\label{eq: betaid combined}
\sum_{e\in \cF_{in}(E)}\int_{e} \beta \cdot n_E = - \sum_{e\in \cF_{out}(E)}\int_{e} \beta \cdot n_E, \quad \text{and} \quad 
\sum_{e\in \cF_{in}(E)}\int_{e} \abs{\beta \cdot n_E} = \sum_{e\in \cF_{out}(E)} \int_{e} \beta \cdot n_E. 
\end{equation}

\begin{definition}[Upwind and downwind traces]\label{def: upwind traces}
    Given $v \in \VhStern$ and an interior face $e\in\cF_i$ with adjacent cells $E_1, E_2$ we follow the usual definitions and define the upwind trace
    \begin{align}\label{eq: upwind/downwind trace}
        v^\uparrow = \begin{cases}
            \trace{v|_{E_1}} & : \quad\beta \cdot n_{E_1} > 0,\\
            \trace{v|_{E_2}} & : \quad\beta \cdot n_{E_2} \ge 0,
        \end{cases}
        \quad \text{and the downwind trace}\quad
        v^\downarrow = \begin{cases}
           \trace{ v|_{E_1}} & : \quad\beta \cdot n_{E_1} < 0,\\
         \trace{   v|_{E_2}} & : \quad\beta \cdot n_{E_2} \le 0,
        \end{cases}
    \end{align}
    where $\trace\cdot$ denotes the trace operator $H^1(E_i) \rightarrow L^2(e)$, $i=1,2$, and
    $n_{E_1}$ and $n_{E_2}$ denote the outward pointing unit vectors of $E_1$ and $E_2$ on $e$.
      These notions are well-defined since we assume that $\beta \cdot n$ does not change sign on any face.
    We extend this definition to a boundary cell $E$ with boundary face $e \in \cF_b$ and unit outward normal $n_E$ as follows
    \begin{align}\label{eq: upwind/downwind trace bdy}
        v^\uparrow = \begin{cases}
            \trace{ v|_{E}} & : \quad \beta \cdot n_{E} > 0,\\
            0 & : \quad \beta \cdot n_{E} \le 0,
        \end{cases}
\qquad\text{and}\qquad
        v^\downarrow = \begin{cases}
            0 & : \quad \beta \cdot n_{E} \ge 0,\\
            \trace{ v|_{E}} & : \quad \beta \cdot n_{E} < 0.
        \end{cases}
    \end{align}
    We further introduce the upwind and downwind normal vectors on any interior face $e\in\cF_i$ as
    \begin{align}
        n^\uparrow = \begin{cases}
            n_{E_1} & : \quad\beta \cdot n_{E_1} > 0,\\
            n_{E_2} & : \quad\beta \cdot n_{E_2} \ge 0,
        \end{cases}
        \qquad\text{and}\qquad
        n^\downarrow = \begin{cases}
            n_{E_1} & : \quad\beta \cdot n_{E_1} < 0,\\
            n_{E_2} & : \quad\beta \cdot n_{E_2} \le 0.
        \end{cases}
    \end{align}
    Note that it holds that $\beta \cdot n^\uparrow \ge 0$ and $\beta \cdot n^\downarrow \le 0$.
\end{definition}

\begin{remark}[In/outflow boundary conditions and stabilized cut cells] 
Based on the setup of the ramp test, the cut cells at the starting point and end point of the ramp will always have a diameter of at least length $h$ and will therefore not need stabilization. Given the large amount of technicalities already present in this work, we will focus on the situation where for $E \in \cI$ there holds $\abs{\partial E \cap \partial \Omega} = 0$, or in other words $e_{in},e_{out} \notin \cF_{b}.$ 
\end{remark}

\begin{definition}[$\beta$-weighted mean traces]\label{def: beta-weighted means}
    We define the $\beta$-weighted means of traces as follows:
    %
    Given $v \in \VhStern$, an element $E$, and a face $e\in \partial E$, we define $\betamean{v|_E}{e} \in \mathbb{R}$ as
    \begin{equation}\label{eq:beta-weighted mean traces}
      \betamean{v|_E}{e}
      \int_e |\beta \cdot n| ds := \int_e |\beta \cdot n| \: \trace{v|_E} ds.
    \end{equation}
    Similarly we define $\beta$-weighted means of functions defined on the interface $e$. For example, the $\beta$-weighted mean of the upwind trace $\betamean{v^\uparrow}{e} \in \mathbb{R}$ is defined as
    \begin{equation}\label{eq:beta-weighted mean upwind}
      \betamean{v^\uparrow}{e}
      \int_e |\beta \cdot n| ds := \int_e |\beta \cdot n| \: v^\uparrow ds,
    \end{equation}
    where the upwind trace $v^\uparrow$ is given by \eqref{eq: upwind/downwind trace} and \eqref{eq: upwind/downwind trace bdy}.
\end{definition}

\begin{remark}
Note that the quantities $\betamean{v|_E}{e}$ and $\betamean{v^\uparrow}{e}$ are constants on each face $e$, whereas the space $\VhStern$ contains elements from $H^1(\Omega)$.
\end{remark}

\begin{definition}[Jump and Average]\label{def: jump and average}
We define the average of $v \in \VhStern$ on face $e= E_1 \cap E_2$ by
\begin{equation*}
    \left\lbrace\mskip-5mu\lbrace{v}\right\rbrace\mskip-5mu\rbrace=\frac 1 2({v}\vert_{_{E_1}}+{v}\vert_{_{E_2}}),
  \end{equation*}
  and the \emph{jump} to be vector-valued given by 
  \begin{equation*}
    \left\llbracket{v}\right\rrbracket:={v}\vert_{_{E_1}}n_{E_1}+{v}\vert_{_{E_2}}n_{E_2},
  \end{equation*}
  with $n_{E_i}$ denoting the outer unit normal vector of cell $E_i$, $i=1,2$. 
  On faces of the physical boundary of the domain $\Omega$, we simply define 
  $ \left\lbrace\mskip-5mu\lbrace{v}\right\rbrace\mskip-5mu\rbrace = v \vert_{_{E}}$ and $\jump{v} = v \vert_{_{E}}n_{E}$.

  We further define jumps and averages of $\beta$-means of traces as follows
  \begin{align*}
    \jump{\betamean{v}{e}} &= \betamean{{v}\vert_{_{E_1}}}{e}n_{E_1}+\betamean{{v}\vert_{_{E_2}}}{e}n_{E_2}
    \quad \text{and}\\
    \average{\betamean{v}{e}} &= \frac 1 2(\betamean{{v}\vert_{_{E_1}}}{e}+\betamean{{v}\vert_{_{E_2}}}{e})
    .
  \end{align*}
\end{definition}
We will later on also use that on interior faces
  \[
  \jump{v}^2 = (({v}\vert_{_{E_1}}-{v}\vert_{_{E_2}})n_{E_1}) \cdot (({v}\vert_{_{E_1}}-{v}\vert_{_{E_2}})n_{E_1}) =
  ({v}\vert_{_{E_1}}-{v}\vert_{_{E_2}})^2 (n_{E_1} \cdot n_{E_1}) = ({v}\vert_{_{E_1}}-{v}\vert_{_{E_2}})^2.
  \]

\begin{notation}[Abbreviation of $\beta$-mean]\label{notation beta-mean}
    Throughout the paper we will \textit{almost always} use $\beta$-weighted means of traces, in particular in our discretization. 
    We introduce the following abbreviated notation, where the element $E$ and face $e$ used in defining the trace are not explicitly stated anymore; 
    their meaning should be clear form the context; the fact that $\beta$-weighted means of traces are used is indicated by an overline notation, i.e.,
    \begin{align*}
        \vweighted &\longrightarrow \bar v
        &\text{and}&&
        \betamean{v^\uparrow}{e} &\longrightarrow \bar v^\uparrow.
    \end{align*}
    {For a stabilized cut cell $E \in \cI$,}
    we reference the $\beta$-weighted mean to particular neighboring cells
    \begin{align*}
        \betamean{v_{E_{in}}}{e_{in}} &\longrightarrow \bar v_{in}
        &\text{and}&&
        \betamean{v|_{E}}{e_{out}}    &\longrightarrow \bar v_{out}.
    \end{align*}
\end{notation}

\subsection{Spatial discretization}\label{subsec: spatial discr}

We start from the DoD-stabilized spatial discretization as proposed in \cite{DoD_2d_linadv_2020}. We replace traces of trial and test functions by their $\beta$-weighted means. As we use piecewise constant discrete spaces, the usage of $\beta$-weighted means does not alter the fully discrete formulation. It changes the stabilization term for the exact solution though. This will be beneficial when we discuss boundedness and consistency later.

The semi-discrete scheme is then given as: Find $u_h(t) \in \mathcal{V}^0$ for all $t \in (0,T)$
such that
\begin{equation}\label{eq: scheme 2d with stab}
        \left(d_t u_h(t),{w_h}\right)_{L^2} + a_h^{upw}\!\left(u_h(t), w_h\right) + J_h\!\left(u_h(t), w_h\right)+ l_h\left(t,w_h\right) = 0
        \qquad \forall ~ w_h \in \mathcal{V}^0,
\end{equation}
with operators $a_h^{upw} : \VhStern\times\mathcal{V}^0 \rightarrow \mathbb{R}$, $J^0_h : \VhStern\times\mathcal{V}^0 \rightarrow \mathbb{R}$, and $l_h : (0,T)\times \mathcal{V}^0 \rightarrow \mathbb{R}$.
The bilinear form $a_h^{upw}$ and the linear form $l_h$ are given by
\begin{align}\label{eq: ah lin adv deg p}
  &\begin{aligned}
      a_h^{upw}(v, w_h):=& \sum_{e \in \cF_b} \int_{e} \left\langle{\beta},{n}\right\rangle^{\oplus}\bar v\, \bar  w_h \text{d}{s}\\
      & +\sum_{e \in \cF_i}\int_{e} \left(\left\langle\jump{\bar{v}},\average{\bar  w_h}\right\rangle
    + \frac{1}{2}\left\vert{\left\langle{\beta},{n_e}\right\rangle}\right\vert
    \left\langle\jump{\bar v},\jump{\bar  w_h}\right\rangle \right) \text{d}{s}\\
    =&
    \sum_{e \in \cF} \int_e \bar v^\uparrow \beta \cdot \jump{\bar  w_h} \text{d}{s},
    \end{aligned}\\
    &\begin{aligned}\label{eq: ah lin adv deg p l_h}
      l_h(t,w_h):=&-\sum_{e \in \cF_b} \int_{e}\left\langle{\beta},{n}\right\rangle^\ominus g(t) \: \bar  w_h\text{d}{s},
    \end{aligned}
\end{align}
where $\cdot^\oplus$ and $\cdot^\ominus$ denote the positive and negative part of a function, respectively, and are zero otherwise.

\begin{remark}
    For any pair of functions $(v, w_h) \in \VhStern  \times \mathcal V^0$ the bilinear form $a_h^{upw}(v, w_h)$ yields exactly the same result as without the $\beta$-weighted means.
    Note further that $\bar v$ and $\bar w_h$ are constant scalar values, which means that they can be pulled out of any integral.
    Their use is mainly motivated by technical reasons, e.g. the definition of the stabilization term $J_h$ becomes easier, in particular it removes integrals over extrapolated functions of the inflow neighbor on the outflow face $e_{out}$. It further simplifies the manipulations in lemma \refwithname{lem:(3.33)}.
\end{remark}

The stabilization term $J_h: \VhStern \times \mathcal{V}^0 \rightarrow \mathbb{R}$ is given by
\begin{equation}\label{def:stab}
    J_h(v,w_h) = J_h^0(v,w_h)= \sum_{E \in \cI} J_h^{0,E}(v,w_h), 
\end{equation}
with
\begin{align}\label{eq: 2d lin adv J0 deg p}
    J_h^{0,E}(v,w_h) &= \eta_E \int_{e_{out}}
    (\bar v_{in}-\bar v_E) \left\langle{\beta},\jump{\bar w_h}\right\rangle\text{d} s.
\end{align}
Here, $\cI$ is the set of cut cells that we want to stabilize
and $\eta_E$ is a penalty parameter. For piecewise constant polynomials we choose
$\eta_E = 1 - \alpha_E$, where $\alpha_E$ denotes the \textit{capacity} of a cut cell $E$, defined below.

This formulation already takes into account that any cut cell $E$ requiring stabilization has a unique inflow face $e_{in}$ with upwind neighbor $E_{in}$ and outflow face $e_{out}$.
The $\beta$-weighted mean $\bar v_{in}$ is defined as indicated in notation \ref{notation beta-mean} and
definition \ref{def: beta-weighted means} above.
The constant value $\bar v_{in}$ is then employed in the face integral over $e_{out}$. Therefore, the
$\beta$-mean $\bar v_{in}$ takes a similar role as the extrapolation operator $\mathcal{L}^\text{ext}$ in \cite{DoD_2d_linadv_2020}. For discrete functions $v_h$ the operator is equivalent to the choice in \cite{DoD_2d_linadv_2020} and only differs for the exact solution.

\begin{definition}[Capacity]\label{eq: def alpha proof}
We define the capacity $\alpha_E$ of a cut cell $E$ as
\begin{equation}\label{eq:newdefalpha}
\alpha_E =\min \left( \frac{|E|}{\tau h \int_{e_{in}} |\beta \cdot n| } , 1 \right).
\end{equation}
Here, $\tau\in \mathbb{R}_+$ is a scalar factor.
%

Note that this choice for the capacity differs from the definition in \cite{DoD_2d_linadv_2020}; the reason is to avoid a dependence of $\alpha_E$ on the time step length $\Delta t$.
\end{definition}

%

We introduce the combined spatial bilinear form
$\ahDoD: \VhStern \times \mathcal{V}^0 \rightarrow \mathbb{R}$ given by
\begin{align}\label{eq: def DoD upw}
    \ahDoD(v,w_h) 
    :&= a_h^{upw}(v,w_h)  + J_h(v,w_h) \nonumber\\
    &=\sum_{e \in \cF} \int_e \bar v^\uparrow \beta \cdot \jump{\bar  w_h} 
    + (1-\alpha_E) \sum_{E \in \cI} \int_{e_{out}(E)} (\bar v_{in}-\bar v_E) \beta \cdot \jump{\bar w_h}  \nonumber\\
    &=\sum_{e \in \cFwoo} \int_e 
    \bar v^\uparrow \beta \cdot \jump{\bar  w_h} + \sum_{E \in \cI} \int_{e_{out}(E)} \left( \alpha_E \bar v^\uparrow + (1-\alpha_E) \bar v_{in} \right) \beta \cdot \jump{\bar  w_h}.
\end{align}
In the last manipulation we used that $\bar v_E$ corresponds to $\bar v^{\uparrow}$ on $e_{out}$.
%
This
induces the operator 
$\AhDoD : \VhStern \rightarrow \mathcal{V}^0$ defined by
\begin{equation}\label{eq: operator AhDoD}
(\AhDoD  v, w_h)_{L^2} := \ahDoD(v, w_h) \quad \forall w_h \in  \mathcal{V}^0.
\end{equation}


\subsection{Definition of norms}

Throughout this work, we will use the following seminorm and norms that have been adjusted to the cut cells present for the ramp model mesh, see figure \ref{fig: ramp 2d cut cell}.



\begin{definition}[$\beta$-seminorm]\label{def: seminorm beta 2d P0} We define for $v \in \VhStern$
\begin{align}\label{eq: seminorm beta 2d P0}
\begin{split}
    \seminormbeta{v}^2 :=& \sum_{e \in \cFwoio} \int_e \abs{\beta \cdot n} \jump{\bar v}^2\\
    &+ \sum_{E \in \cI} \left( 
    \alpha_E \sum_{e \in \{e_{in}(E),e_{out}(E)\}}\int_{e} \abs{\beta \cdot n} \jump{\bar v}^2
    + (1-\alpha_E) \int_{e_{out}(E)} \vspace*{-5mm}(\beta \cdot n_{out}) (\bar v_{out}-\bar v_{in})^2  \right).
    \end{split}
    \end{align}
\end{definition}
    Note that the norm uses squares of $\beta$-weighted means and not $\beta$-weighted means of squares. Boundary integrals over faces of stabilized cut cells that involve normal jumps are scaled with capacity $\alpha_E$. The seminorm also contains an \textit{extended} jump, which measures the difference between the traces from the two neighboring cells of a stabilized cut cell, $\bar v_{in}$ and $\bar v_{out}$.

\begin{definition}\label{def: uwbsterm 2D P0}
For $v \in \VhStern$ we define
\[
\normuwb{v}^2 := \norm{v}_{L^2(\Omega)}^2 + \seminormbeta{v}^2
\]
and
\begin{equation}
    \normuwbStern{v}^2:= \normuwb{v}^2 + \sum_{E \not\in \cI} \|\, |\beta \cdot n |^{1/2} \bar v \|_{L^2(\partial E)}^2
    + \sum_{E \in \cI} \alpha_E \|\, |\beta \cdot n |^{1/2} \bar v\|_{L^2(\partial E)}^2.
\end{equation}
\end{definition}

\subsection{Fully discrete scheme and main result}\label{subsec: properties bilin}
We combine the spatial discretization described in section \ref{subsec: spatial discr} with explicit Euler in time. 
Starting from the given initial data $u_0$, we compute discrete initial data $u_h^0 \in \mathcal{V}^0$ as
$u_h^0 = \Pi_h u_0$ where $\Pi_h: L^2(\Omega) \to \mathcal{V}^0$ denotes the $L^2$ projection operator.
Denoting by $t^1,\ldots,t^N$ the time steps, this leads to the following fully discrete scheme:
Given $u_h^0 \in \mathcal{V}^0$, find $\{u_h^n\}_n \in \mathcal{V}^0$ such that
\begin{equation}\label{eq: scheme 2d with stab fully discrete}
        \left(\frac{u_h^{n+1}- u_h^n}{\Delta t},{w_h}\right)_{L^2} + \ahDoD\!\left(u_h^n, w_h\right) + l_h\left(t^n,w_h\right) = 0
    \end{equation}
for all $w_h \in \mathcal{V}^0 $ and for all $n=0,\ldots,N-1.$


\paragraph*{Main result {\normalfont (c.f. theorem \ref{thm:main})}:} We will show a result of the form: for all $M$ with $M \dt\le T$ there holds for any fixed $0<\varepsilon < \frac 1 2$ under a constraint on the time step length of the form $\dt \le \frac{1-2\varepsilon}{(1+\varepsilon)\max(4\norm{\beta}_{\infty},\tau^{-1})} h$
 \begin{equation*}
  \| u(t^M,\cdot) - u_h^M \|_{L^2(\Omega)}^2 + \sum_{n=0}^{M-1} \varepsilon \dt \seminormbeta{u(t^n,\cdot) - u_h^n}^2
  \leq (1+\epsinv)[Ch \norm{u}^2_{C^0([0,T],H^1(\Omega))} + C\dt^2 \norm{u}^2_{C^2([0,T],L^2(\Omega))} ].
 \end{equation*}
 We stress that the CFL condition is independent of the size of small cut cells. Further, as verified in the numerical results in section \ref{sec: numerical results}, the estimate is optimal with respect to the error measured in $\seminormbeta{\cdot}$.

The proof of this result and more details about the constant $C$ and the prerequisites are given in section~\ref{sec: proof of main result}.

To derive this estimate we 
build on the powerful framework presented in \cite{BURERNFER2010}. The proof exploits several properties of the bilinear form $\ahDoD$, which we will prove in section \ref{sec: auxiliary results} and are summarized here:
\begin{description}
    \item[\refastitle{lemma: consistency_new}:]
    \begin{equation*}
        \abs{J_h(v, w_h)} \leq \sqrt{\tau  h} \|\beta \|_{W^{1,\infty}} \|v\|_{H^{1}(\Omega)} \seminormbeta{w_h}
        \qquad
        \forall ~ (v, w_h) \in H^1(\Omega) \times \mathcal{V}^0.
    \end{equation*}
    %
    \item[\refastitle{lemma: coercivity 2d lin adv P0}:]
    \begin{equation*}
        a_h^{DoD}(v_h, v_h) = \frac{1}{2}|v_h|^2_{\beta} \geq 0
        \qquad \forall~ v_h \in \mathcal{V}^0.
    \end{equation*}
 %
    \item[\refastitle{lem:(3.33)}:]
    \begin{equation*}
    \abs{a_h^{DoD}(v, w_h)} \leq \normuwbStern{v} |w_h|_{\beta}
    \qquad
    \forall ~ (v, w_h) \in \mathcal{V}_{*}^0 \times \mathcal{V}^0.
    \end{equation*}
    \item[\refastitle{lem:1bound}:]
\begin{equation*}
    \| \AhDoD v_h \|_{L^2(\Omega)}  \leq 
    \sqrt{\frac{\Ctr}{h}}\seminormbeta{v_h}
    \qquad
    \forall ~ v_h \in \mathcal{V}^0,
\end{equation*}
with
$A_h^{DoD}$ as introduced in \eqref{eq: operator AhDoD},
and $\Ctr$ independent of $h$, $\Delta t$, $\alpha_E$, and the size of small cut cells (but depending on $\beta$ and $\tau$).
\end{description}

\begin{notation}[Constants] Throughout this work, we will use both named constants, e.g., $\Ctr$ in lemma \ref{lem: inv est const}, as well as a generic constant $C$, whose value may change from line to line in estimates. The constants are generally independent of $h, \dt,\alpha_E$, the size of small cut cells, and the exact solution $u$. Whether they depend on other quantities, e.g. $\beta$ or $\tau$, depends on the context.
\end{notation}

\section{Geometric cut cell specific estimates}\label{sec:geom}


To prove the auxiliary results in section \ref{sec: auxiliary results}, one usually needs several standard results, in particular inverse estimates and estimates on the projection error. We cannot use standard variants but need versions that account for the fact that the cut cell mesh is not quasi-uniform, and are adapted to the roles of $|\beta \cdot n|$ and $\alpha_E$ in the bilinear form. In particular, we will exploit that boundary integrals typically contain the integrand $\abs{\beta \cdot n}$. These results are independent of our particular stabilization.

\begin{lemma}[Inverse trace estimate]\label{lem: inv est const}
Let the set of stabilized cut cells $\cI$ be given by \eqref{eq: set of stab cut cells}.
For a non-stabilized cell $E$, including non-stabilized cut cells, there holds
\begin{equation}\label{eq:inv est const unstab}
\sum_{e \in \cF_{in}(E)} \int_e \abs{\beta(x) \cdot n(x)} = \sum_{e \in \cF_{out}(E)} \int_e \abs{\beta(x) \cdot n(x)} \le \frac{\Ctr[1]}{h} \abs{E}
 \qquad \text{with } \Ctr[1] = 4 \norm{\beta}_{\infty}.
\end{equation}
For a stabilized cut cell $E$, there holds with $\alpha_E$ and $\tau$ from \eqref{eq:newdefalpha} 
\begin{equation}\label{eq:inv est const stab}
 \alpha_E \sum_{e \in \cF_{in}(E)} \int_e |\beta(x) \cdot n(x)| = \alpha_E \sum_{e \in \cF_{out}(E)} \int_e |\beta(x) \cdot n(x)| \le \frac{\Ctr[2]}{h} \abs{E}
 \qquad \text{with } \Ctr[2] = \frac{1}{\tau}.
\end{equation}
For later use of combining the two estimates we define the constant $\Ctr := \max(\Ctr[1], \Ctr[2])$.
\end{lemma}

\begin{proof}
We generally have four kinds of cells: Cartesian background cells, five-sided cut cells, four-sided cut cells, and three-sided cut cells. 
\begin{itemize}
\item 
Cartesian background cells: for a Cartesian cell $E$, we have by using \eqref{eq: def beta-inf}
\[   \sum_{e \in \cF_{in}(E)} \int_e \abs{\beta(x) \cdot n(x)} \le \sum_{e \in \cF_{in}(E)}\int_e \abs{n(x)}_2\abs{\beta(x)}_{2} 
\le |\cF_{in}(E)| \norm{\beta}_{\infty} h, \]
as well as
\[   \sum_{e \in \cF_{out}(E)} \int_e \abs{\beta(x) \cdot n(x)} \le \sum_{e \in \cF_{out}(E)}\int_e \abs{n(x)}_2\abs{\beta(x)}_{2} 
\le |\cF_{out}(E)| \norm{\beta}_{\infty} h. \]
Here, $\abs{\cF}$ denotes the cardinality of the set $\cF$.
Due to \eqref{eq: betaid combined} we can combine these estimates to get
\[   \sum_{e \in \cF_{in}(E)} \int_e \abs{\beta(x) \cdot n(x)} = \sum_{e \in \cF_{out}(E)} \int_e \abs{\beta(x) \cdot n(x)}
\le \min(|\cF_{in}(E)|, |\cF_{out}(E)|) \norm{\beta}_{\infty} h    \le  \frac{2 \norm{\beta}_{\infty}}{h} \abs{E}. \]
In the last step, we used $\min(|\cF_{in}(E)|, |\cF_{out}(E)|)\leq 2$ and $\abs{E} = h^2$.

\item Five-sided cut cells: for a five-sided cut cell $E$, $\abs{E} > \frac 1 2 h^2$. Together with the fact that the cut face is a no-flow face, this implies that we have the same situation as for Cartesian cells except for losing a factor of 2.
\item Four-sided cut cells: on the cut boundary of a four-sided cut cell
it holds that
$\beta \cdot n = 0$. We denote by $e_3$ the face opposite to the cut boundary, which has $\abs{e_3} = h$, and by $e_1$ and $e_2$ the two intersected faces for which we have $0 < \abs{e_1},\abs{e_2}<h$. There holds $\abs{E} = \frac{1}{2}(\abs{e_1}+\abs{e_2})h.$ We focus on the face $e_3$ as this includes the other two faces. By incompressibility
\[
\int_{e_3} \abs{\beta \cdot n} \le \int_{e_1} \abs{\beta \cdot n} + \int_{e_2} \abs{\beta \cdot n}
\le \norm{\beta}_{\infty} (\abs{e_1}+\abs{e_2})  = 
2\frac{\norm{\beta}_{\infty}\abs{E}}{h}.
\]
Therefore,
\[
\sum_{e \in \cF_{in}(E)} \int_e \abs{\beta(x) \cdot n(x)} = \sum_{e \in \cF_{out}(E)} \int_e \abs{\beta(x) \cdot n(x)} 
\le \sum_{i=1}^4 \int_{e_i} \abs{\beta(x) \cdot n(x)}
\le 4\frac{\norm{\beta}_{\infty}\abs{E}}{h}.
\]
\item Three-sided cut cells: These are always right-angled triangles. The hypotenuse is formed by the cut and thus is a boundary face with $\beta \cdot n=0$. As a result, the sets $\cF_{in}(E)$ and $\cF_{out}(E)$ consist of only one face each, which are the two catheti of the triangle. Let us call these two faces $e_1$ and $e_2$. We will now prove the claimed estimates for $e_1$ and $e_2$. It holds that $\abs{E} = \frac{1}{2}\abs{e_1}\abs{e_2}.$ W.l.o.g. we assume $\abs{e_1} \le \abs{e_2}$. By means of the incompressibilty, we can rewrite
\[
\int_{e_2} \abs{\beta \cdot n} = \int_{e_1} \abs{\beta \cdot n} \le \abs{e_1} \norm{\beta}_{\infty}  = 2\frac{\norm{\beta}_{\infty}}{\abs{e_2}} \abs{E}.
\]
Note that this implies that all considerations are also true for the shorter face $e_1$. 
Looking now at the situation of a non-stabilized cell, we know that $\abs{e_2} \ge \frac 1 2 h$ and therefore we get
\[
\int_{e_2} \abs{\beta \cdot n} \le 4\frac{\norm{\beta}_{\infty}}{h} \abs{E}.
\]
On a stabilized cut cell we have by definition of $\alpha_E$, compare \eqref{eq:newdefalpha}, 
\begin{equation}
\alpha_E =\min \left( \frac{|E|}{\tau h \int_{e_{in}} |\beta \cdot n| } , 1 \right) \le
\frac{|E|}{\tau h \int_{e_{in}} |\beta \cdot n| }.
\end{equation}
This together with the incompressibility condition imply for both $i=1$ and $i=2$
\[
\alpha_E \int_{e_i} \abs{\beta \cdot n} \le \frac{|E|}{\tau h \int_{e_{in}} |\beta \cdot n| } \int_{e_i} \abs{\beta \cdot n}
= \frac{\abs{E}}{\tau h}.
\]
\end{itemize}
\end{proof}

\begin{lemma}[Inverse estimate]\label{lemma: inverse eq betaseminorm}
There holds for $w_h \in \mathcal{V}^0$
\begin{equation}\label{eq:inverse seminormbeta}
    \seminormbeta{w_h} \leq
    2\sqrt{\frac{\Ctr}{h}} \|w_h\|_{L^2},
\end{equation}
with 
$\Ctr$ from lemma \ref{lem: inv est const}.
\end{lemma}

\begin{proof}
By definition
\begin{displaymath}
    \seminormbeta{w_h}^2 = T_1 + T_2 + T_3
\end{displaymath}
with
\begin{align*}
    T_1 & = \sum_{e \in \cFwoio} \int_e \abs{\beta \cdot n} \jump{\bar w_h}^2, \qquad T_2 =  \sum_{E \in \cI} \left( \alpha_E \int_{e_{out}(E)} \abs{\beta \cdot n} \jump{\bar  w_h}^2 + \alpha_E \int_{e_{in}(E)} \abs{\beta \cdot n} \jump{\bar  w_h}^2 \right),\\
    T_3 & = \sum_{E \in \cI} \left( (1-\alpha_E) \int_{e_{out}(E)} (\beta \cdot n_{out}) (\bar w_{out}-\bar w_{in})^2  \right).
\end{align*}
By $(\bar w_{E_1} - \bar w_{E_2})^2 \le 2(\bar w_{E_1})^2 + 2(\bar w_{E_2})^2$ and by using the divergence theorem to rewrite integrals over $e_{out}$ as integrals over $e_{in}$, we get
\[
T_3
\le 2 \sum_{E \in \cI} \left( (1-\alpha_E) \int_{e_{out}} (\beta \cdot n_{out}) (\bar w_{out})^2 +
(1-\alpha_E) \int_{e_{in}} \abs{\beta \cdot n_{in}} (\bar w_{in})^2 \right).
\]
Note that all involved $\beta$-means $\bar w_h$ (replacing traces of $w_h$) in $T_3$ belong to unstabilized cut cells $E \not \in \cI$.
Using $\jump{\bar  w_h}^2 = (\bar  w_{E_1} - \bar  w_{E_2})^2 \le 2(\bar  w_{E_1})^2 + 2(\bar  w_{E_2})^2$ again, we get for $T_2$
\[
T_2 \le 2 \sum_{E \in \cI} \left( \alpha_E \int_{e_{out}} (\beta \cdot n_{out}) ((\bar w_E)^2+ (\bar w_{out})^2) +
\alpha_E \int_{e_{in}} \abs{\beta \cdot n_{in}} ((\bar w_E)^2 + (\bar w_{in})^2) \right).
\]
Therefore, using that $\beta \cdot n = 0$ on the diagonal face of a stabilized cut cell, there holds
\[
T_2 + T_3 \le 2 \sum_{E \in \cI} \left( \alpha_E \int_{\partial E} \abs{\beta \cdot n} (\bar w_E)^2   \right)
+ 2 \sum_{E \in \cI} \left(  \int_{e_{out}} (\beta \cdot n_{out}) (\bar w_{out})^2 +
 \int_{e_{in}} \abs{\beta \cdot n_{in}} (\bar w_{in})^2 \right).
\]
For $T_1$, the sum over boundary integrals of Cartesian cells and non-stabilized cut cells, we again employ that
$\jump{\bar w_h}^2 \le 2(\bar w_{E_1})^2 + 2(\bar w_{E_2})^2$. Further note that the missing boundary integrals are exactly the boundary integrals given in $T_2+T_3$. Therefore,
\[
T_1 + T_2 +T_3 \le
2 \sum_{E \not \in \cI} 
\int_{\partial E} |\beta \cdot n|(\bar w_E)^2 + 
2 \sum_{E \in \cI} 
\alpha_E \int_{\partial E} |\beta \cdot n| (\bar w_E)^2.
\]
We now exploit that $w_h \in \mathcal{V}^0$ is piecewise constant and that therefore the estimates in lemma \ref{lem: inv est const} hold in the same way for boundary and volume integrals over $w_h$. Note that
$\int_{\partial E} |\beta \cdot n| = \sum_{e \in {\cF_{in}(E)}\cup \cF_{out}(E)} \int_e \abs{\beta \cdot n}$.
By lemma \ref{lem: inv est const}, the observation that $|E|\bar w_E^{2} =\|w_h\|_{L^2(E)}^{2}$, and the definition of $\Ctr := \max(\Ctr[1], \Ctr[2])$ made earlier 
we conclude that
\begin{equation}\label{eq: last step pf inv est}
\seminormbeta{w_h}^2 =  T_1 + T_2 + T_3 
\le 4 \frac{\Ctr[1]}{h} \sum_{E \not \in \cI} |E|\bar w_E^{2} + 4 \frac{\Ctr[2]}{h} \sum_{E\in \cI} |E|\bar w_E^{2}
\le 4 \frac{\Ctr}{h} \sum_{E \in \cT_h} \| w_h \|_{L^2(E)}^2
\end{equation}
to finish the proof.
\end{proof}

\subsection{Controlling projection errors}

Due to the presence of the small cut cells, the standard projection error estimates do not apply here. We therefore provide the following results.

\begin{lemma}[Projection Error]\label{lem:est_xipi}
Let us assume that there is  a (uniform) constant $c_b>0$ such that for each inflow and each outflow face of a stabilized cut cell $E \in \cI$ there holds $|\beta \cdot n| > c_b$ pointwise. Denote by $\xi_\pi(t,\cdot) = u(t,\cdot) - \Pi_h u(t,\cdot)$ the error of the $L^2$ projection operator $\Pi_h: H^1(\Omega) \to \mathcal{V}^0$. 
Then there holds for $u \in C^0([0,T],H^1(\Omega))$ 
\[ \| \xi_\pi(t,\cdot)  \|_{L^2(\Omega)} \leq  \frac{\sqrt{2}}{\pi} h\| \nabla u(t,\cdot) \|_{L^2(\Omega)},  \]
and there exists a generic constant $C$, depending on $\tau$, $\beta$, and the ramp angle $\gamma$, but independent of $h$, the mesh refinement process, and the size of small cut cells, such that
\[ \normuwbStern{ \xi_\pi(t,\cdot)  } \leq  C \sqrt{h} \| \nabla u(t,\cdot) \|_{L^2(\Omega)}.
\]
\end{lemma}

\begin{remark}
The technical assumption in lemma \ref{lem:est_xipi} is satisfied if, e.g., the mesh size is fine enough, the ramp angle $\gamma$ is bounded away from $0^{\circ}$ and $90^{\circ}$, and the velocity field $\beta$ does not vanish along the ramp. 
\end{remark}

\begin{proof}
    By construction, the mean of $\xi_\pi$ is zero on each cell. This allows to use cell-wise Poincare inequalities. Since each cell is convex we may use the optimal Poincare constant for convex domains derived in \cite{Payne_1960}. This implies 
    \begin{equation}\label{eq:loc_approx}\| \xi_\pi(t,\cdot)  \|_{L^2(E)} \leq  \frac{\text{diam}(E)}{\pi}\| \nabla u(t,\cdot) \|_{L^2(E)} \quad
    \forall E \in \cT_h.
    \end{equation}
  Squaring, summing this estimate over all cells, and bounding the local cell size $\text{diam}(E)$ by the global cell size $\sqrt{2}h$, with $h$ being the face length of a Cartesian background cell, implies the first assertion. 

For proving the second estimate we observe that
$\xi_\pi(t,\cdot)\in \VhStern$. Thus we consider functions $v\in \VhStern$.
Recalling the definition of the norm $\normuwbStern{\cdot}$ in definition \ref{def: uwbsterm 2D P0}
\begin{equation}
    \normuwbStern{v}^2:= \underbrace{\norm{v}_{L^2(\Omega)}^2}_{T_1} + \underbrace{\seminormbeta{v}^2}_{T_2} + \underbrace{\sum_{E \not\in \cI} \|\, |\beta \cdot n |^{1/2} \bar v \|_{L^2(\partial E)}^2}_{T_3}
    + \underbrace{\sum_{E \in \cI} \alpha_E \|\, |\beta \cdot n |^{1/2} \bar v\|_{L^2(\partial E)}^2}_{T_4},
\end{equation}
we first observe that for $v = \xi_\pi(t,\cdot)$ the term $T_1$ is bounded according to \eqref{eq:loc_approx}.
The $\beta$-seminorm term $T_2$ can be bounded, up to a constant $C$, by an expression of the form of $T_3 + T_4$ 
by replacing squares of jumps $\jump{\cdot}^2$ by sums of squares. 
For example, considering the first term in $\seminormbeta{\cdot}$, i.e., a face $e \in \cFwoio$, and using as we did before $\jump{\bar v}^2 \le 2(\bar v_{E_1})^2 + 2(\bar v_{E_2})^2$, we have
\begin{equation*}
\int_e \abs{\beta \cdot n} \jump{\bar v}^2 \le \int_e \abs{\beta \cdot n} \left( 2 (\bar v_{E_1})^2 + 2(\bar v_{E_2})^2  \right) \le 2 \int_e \abs{\beta \cdot n}(\bar v_{E_1})^2 + 2 \int_e \abs{\beta \cdot n}(\bar v_{E_2})^2,
\end{equation*}
which correspond, up to the factor of 2, to terms showing up in $T_3$. The $\alpha_E$-weighted terms in $\seminormbeta{\cdot}$ are treated similary, with some terms corresponding to terms in $T_3$ (using $\alpha_E \le 1$) and some to terms in $T_4$. The $(1-\alpha_E)$-weighted terms in $\seminormbeta{\cdot}$ finally, which contain extended jumps, follow the same logic but use in addition that for a stabilized cut cell $E$ by incompressibility, one can replace
$\int_{e_{out}(E)} \abs{\beta \cdot n} \bar v_{in}^2$ by $\int_{e_{in}(E)} \abs{\beta \cdot n} \bar v_{in}^2$.

It is now evident that the main task is to bound $T_3$ and $T_4$. We begin by proving that local estimates for traces immediately lead to bounds for $\beta$-weighted means of traces. We have
\begin{align*}
 0 \le \int_e (v- \overline{v})^2 |\beta \cdot n| ds 
= \int_e (v^2- 2\overline{v}^2 + \overline{v}^2) |\beta \cdot n| ds
= \int_e v^2 |\beta \cdot n| ds
- \int_e  \overline{v}^2 |\beta \cdot n| ds,
\\
\Rightarrow\qquad \int_e \overline{v}^2 |\beta \cdot n| ds \leq \int_e v^2 |\beta \cdot n| ds.
\end{align*}
Therefore, it is sufficient to derive estimates for actual traces of $v$.




We now provide some geometric considerations.
Recall that $0<\gamma< \pi/2$ denotes the ramp angle in radians and let
$E$ be a non-stabilized cut cell. Then $E$ contains a right triangle with length of the longest leg $\ell \geq h/2$
and angle between longest leg and hypothenuse $\delta$ with $\delta= \gamma$ if $\gamma < \pi/4$ and $\delta= \pi/2 - \gamma$ otherwise. This is sketched in figure \ref{fig: geom proj error}. 
Then the length of the faces of $E$ are given by $\ell$, $\ell \tan \delta$, and $\tfrac{\ell}{ \cos \delta } $.

\begin{figure}[b]
  %
  \begin{subfigure}[b]{0.50\textwidth}
  \centering
  \begin{tikzpicture}[scale=1.0]
\draw[] (-0.2,0) -- (3.2,0);
\draw[] (0,0.1) -- (0,-0.8);
\draw[] (2,0.1) -- (2,-0.3);
\draw[dotted] (0,-0.3) -- (2,-0.3);
\draw[thick] (-0.2,-0.85) -- (3.4,0.05);
\draw [<->] (2,0.2) -- (3.2,0.2);
\node at (2.6,0.4){\small $\ge \tfrac h 2$};
\draw [<->] (0,0.2) -- (2,0.2);
\node at (1,0.4){\small $h$};
\fill[MyColorBlue] (-0.2,-0.85) --
      (-0.2,-1.2) --
      (3.4,-1.2) --
      (3.4,0.05) --
      (-0.2,-0.85);
\end{tikzpicture}
  \caption{Four-sided and non-stabilized three-sided cut cell.}
  \label{fig: geom proj error a}
    \end{subfigure}%
\begin{subfigure}[b]{0.50\textwidth}
  \centering
  \begin{tikzpicture}[scale=1.2]
\draw[] (0,0) -- (3,0);
\draw[] (0,0) -- (0,-0.75);
\draw[] (0,-0.75) -- (3,0);
\draw[dashed](0,-0.75) -- (1,-0.75);
\draw[] (1,-0.75) arc (0:14:1cm);
\draw[<-] (0.8,-0.65) -- (1.05,-0.85);
\node[] at (1.2,-0.95) {$\gamma$} ;
\draw[] (2,0) arc (180:194:1cm);
\draw[->] (1.8,-0.125) -- (2.2,-0.125);
\node[] at (1.6,-0.15) {$\delta$} ;
\node[] at (1.4,0.3) {$\ell$} ;
\node[] at (-0.6,-0.3) {$\ell \tan\delta$} ;
\node[] at (1.8,-0.6) {$\tfrac{\ell}{\cos\delta}$} ;
\end{tikzpicture}
  \caption{Information for the contained triangle.}
  \label{fig: geom proj error b}
    \end{subfigure}%
  \caption{Geometric considerations for proving projection errors for ramp angle $\gamma < \pi/4$: All non-stabilized cut cells contain triangular cells with the longer leg having a length $\ell \ge \tfrac{h}{2}$. Shown in (a) are examples for four-sided and three-sided cut cells. Figure (b) shows detailed information for this triangle.}
  \label{fig: geom proj error}
\end{figure}
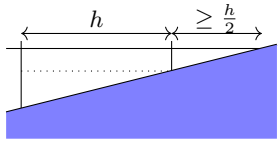
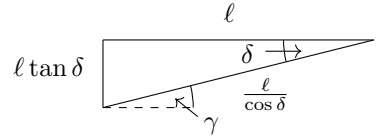

The incircle radius $r$ is given by
\[
r = \sqrt{\frac{\ell^2}{4} 
\frac{(-1+ \tan \delta + \frac{1}{\cos \delta}) (1- \tan \delta + \frac{1}{\cos \delta})(1+ \tan \delta - \frac{1}{\cos \delta})}{1+ \tan \delta + \frac{1}{\cos \delta}}  }
=: \ell c(\delta) \geq \frac{h}{2}c(\delta)
\]
 where $c(\delta)$ is some number that depends solely on $\gamma$, since $\delta$ is determined by $\gamma$. 
    For each non-stabilized cut cell $E$ there exists a point $x_E$ such that
    $( x- x_E) \cdot n \geq \tfrac{h}{2} c(\delta)$ for all $x \in  \partial E$, where $n$ denotes the outward unit normal to $E$. The reason is that for example for the incenter $\hat{x}$ there holds $(x-\hat{x}) \cdot n \geq r$.
    

    We now estimate term {$T_3$:} On these unstabilized cells we define the function $\sigma(x):= (x-x_E)/h$. Note that $\tfrac{2}{ c(\delta)}\sigma(x) \cdot n \ge 1.$
    Then, for each $v \in H^1(E)$ we
    estimate the element-wise $L^2$-boundary-norm as follows
    \begin{multline}
    \|\, |\beta \cdot n |^{1/2} \bar v \|_{L^2(\partial E)}^2 \le
    \|\, |\beta \cdot n |^{1/2} v \|_{L^2(\partial E)}^2 = 
        \int_{\partial E} |v|^2 |\beta \cdot n|\\
        \leq \frac{2}{ c(\delta)} \| \beta\|_\infty \int_{\partial E}  |v|^2 \sigma \cdot n
        = \frac{2}{ c(\delta)} \| \beta\|_\infty \int_E \nabla \cdot( |v|^2 \sigma)dx\\
        \leq \frac{4}{ c(\delta)} \| \beta\|_\infty ( \sqrt{2}\| v\|_{L^2(E)} \| \nabla v\|_{L^2(E)} + h^{-1} \| v\|_{L^2(E)}^2),
    \end{multline}
    where we have used
    \begin{equation}\label{eq:chainrulediv}
     \int_E \nabla \cdot( |v|^2 \sigma) \: dx
      =\int_E ( (\nabla \cdot \sigma) |v|^2 
      + 2 v \sigma \cdot \nabla v  ) \: dx
    \end{equation}
    with  $\| \sigma \|_{L^\infty(E)}  \leq \sqrt{2}$ and $\| \nabla \cdot \sigma \|_{L^\infty(E)} \leq \frac{2}{h}.$
    Therefore, by using \eqref{eq:loc_approx}, we get with $\text{diam}(E)\le \sqrt{2}h$
    \begin{equation}
        \int_{\partial E} | \bar \xi_\pi|^2 |\beta \cdot n|
        \leq \frac{4}{ c(\delta)} \| \beta\|_\infty \left( \frac{2h}{\pi} \| \nabla u(t,\cdot)\|^2_{L^2(E)} + 
        \frac{2h}{\pi^2} \| \nabla u(t,\cdot)\|_{L^2(E)}^2\right)
        \leq \frac{4}{c(\delta)} h  \| \beta\|_\infty \| \nabla u(t,\cdot)\|_{L^2(E)}^2.
      \end{equation}  
      We note that $c(\delta)$ only depends on the ramp angle and therefore does not change under mesh refinement. 
      The case of a Cartesian cell is covered in the argumentation above by replacing $c(\delta)$ with 1.
    This finishes the estimate for $T_3$.

    Finally we estimate the last term {$T_4$:}
    For each stabilized cut cell $E$ we show how to handle the inflow face, the outflow face can be handled analogously. 
    The argument follows \cite[Lem. 1.49]{DiPietro_Ern}. Each such cell is a triangle and we denote 
    \[ \sigma(x)= \frac{|e_{in}|}{2|E|} (x- x_E) (-\beta(x) \cdot n_{in}),
    \]
    where $x_E$ is the vertex of $E$ opposite of $e_{in}$ and $n_{in}$ denotes the normal vector to $e_{in}$ pointing away from $E$ (which is constant).
    Note that $\sigma \cdot n$ is $|\beta \cdot n|$ on $e_{in}$ and zero on the other faces of $E$ (as for the other two faces, $(x-x_E)$ is tangential to the corresponding face and therefore perpendicular to the normal vector of that face).
    Then
    \begin{align*}
    \alpha_E \int_{e_{in}} &|\bar v|^2 |\beta \cdot n|
    \leq \alpha_E \int_{e_{in}} |v|^2 |\beta \cdot n|
    = \alpha_E \int_{\partial E} |v|^2 \sigma \cdot n= \alpha_E \int_E \nabla \cdot (|v|^2 \sigma)
     \\
     &= \alpha_E \int_E ( 2 v \sigma \cdot \nabla v + (\nabla \cdot \sigma) |v|^2 
      ) \: dx 
    \le \alpha_E \left( 2\| v\|_{L^2(E)} \|\nabla v\|_{L^2(E)} \norm{\sigma}_{L^{\infty}(E)} + \| v\|_{L^2(E)}^2\norm{\nabla \cdot\sigma}_{L^{\infty}(E)} \right) \\
    &\le 2 \alpha_E \left(\| v\|_{L^2(E)} \|\nabla v\|_{L^2(E)} \|\beta \|_\infty \frac{h |e_{in}|}{2|E|} + \| v\|_{L^2(E)}^2\left(\frac{|e_{in}|}{2|E|} \norm{\beta}_{\infty} + \frac{h |e_{in}|}{4|E|} \|\nabla \beta\|_{\infty}\right)\right).
    \end{align*}
    Here, we used $\abs{x-x_E}_2 \le \text{diam}(E) \le \sqrt{2} \frac{h}{2} \le h.$
     Note that by definition of $\alpha_E$ 
    \[ \alpha_E\frac{h |e_{in}|}{2|E|} \leq
    \frac{|E|}{\tau h \int_{e_{in}} |\beta \cdot n| ds } \frac{h |e_{in}|}{2|E|}  = 
    \frac{|e_{in}|}{2\tau \int_{e_{in}} |\beta \cdot n| ds }  \leq \frac{1}{2\tau c_b},
    \]
    with the constant $c_b$ from the assumption of the lemma,
    so that 
    \[
    \alpha_E \int_{e_{in}} |\bar v|^2 |\beta \cdot n|
    \le \frac{1}{\tau c_b} \left(\| v\|_{L^2(E)} \|\nabla v\|_{L^2(E)} \norm{\beta}_{\infty}  +  \| v\|_{L^2(E)}^2
    (h^{-1}\norm{\beta}_{\infty} +  \|\nabla \beta\|_{\infty} )\right). \]
 Together with \eqref{eq:loc_approx}, this implies
       \begin{equation}
    \alpha_E \int_{e_{in}} |\bar \xi_\pi(t,\cdot)|^2 |\beta \cdot n|
        \leq \frac{h}{\tau c_b} \| \nabla u(t,\cdot)\|_{L^2(E)}^2 \left(\norm{\beta}_{\infty} + \frac{h}{\pi^2} \norm{\nabla \beta}_{\infty}  \right).   
      \end{equation}   
This finishes the proof of the second assertion. 
\end{proof}

\section{Auxiliary results}\label{sec: auxiliary results}

In this section we prove the central properties of the bilinear form that we summarized in section 
\ref{subsec: properties bilin} above.

\subsection{Consistency}\label{subsec: consistency}


See \eqref{def:stab} and \eqref{eq: 2d lin adv J0 deg p} for the definition of the stabilization term.

\begin{lemma}[Consistency]\label{lemma: consistency_new}
  For each $v \in H^1(\Omega)$ and $w_h \in \mathcal{V}_0$ there holds
\begin{equation}\label{eq: consistency 2}
\abs{J_h(v,w_h)} \le \sqrt{\tau  h} \norm{\beta}_{W^{1,\infty}} \|v\|_{H^{1}(\Omega)} \seminormbeta{w_h}.
\end{equation}
\end{lemma}


\begin{proof}
    Let $E \in \cI$ be an arbitrary stabilized cut cell with $\alpha_E < 1$. 
    (For $\alpha_E=1$, the stabilization term $J_h^{0,E}$ vanishes.)
    Then
    \[ 
    J_h^{0,E}(v,w_h) = (1 - \alpha_E)\int_{e_{out}} (\bar v_{in} - \bar v_E) \beta \cdot \jump{\bar w_h} ds.
    \]
    We define $k_E:= \int_{e_{out}} ( \bar v_{in} - \bar v_E) \beta \cdot n_E \: ds$.
    Then
    $|J_h^{0,E}(v,w_h)|= (1-\alpha_E)\abs{k_E} \abs{\jump{\bar w_h}}_2$
    and
    using that $\nabla \cdot \beta = 0$
    the absolute value of $k_E$ can be bounded as follows
    \begin{multline}\label{eq:kEcontrol}
       \abs{k_E}=  \abs{\int_{e_{out}} ( \bar v_{in} - \bar v_E) \beta \cdot n_E \: ds}
       = \abs{-  \int_{e_{out}}  \bar v_E \beta \cdot n_E \: ds-  \int_{e_{in}}  \bar v_{in} \beta \cdot n_E \: ds}
       \\
       = \abs{-  \int_{e_{out}}  v \beta \cdot n_E \: ds-  \int_{e_{in}}  v \beta \cdot n_E \: ds}
       = \abs{\int_{\partial E} v \beta \cdot n_E \: ds}
       = \abs{ \int_E \nabla \cdot (v \beta)} \leq \|\beta \|_{W^{1,\infty}} \| v \|_{W^{1,1}(E)},
    \end{multline}
    where we have used that $v \in H^1(\Omega)$ ensures that $v$ has a  unique trace on each face. 

    Since $\jump{\bar w_h}$ is constant on each face, we have
    \begin{multline}
      | J_h(v,w_h)|
      \le \sum_{E \in \cI}\left| J_h^{0,E}(v,w_h) \right|\\
      = \sum_{E \in \cI}\left| (1-\alpha_E) \int_{e_{out}} (\bar v_{in} - \bar v_E) \beta \cdot \jump{\bar w_h} \: ds \right|
      = \sum_{E \in \cI} (1-\alpha_E) \abs{k_E}  \abs{ \jump{\bar w_h}}_2 \: \\
       \le {\underbrace{\left( \sum_{E \in \cI} \frac{1}{\alpha_E \int_{e_{out}}{|\beta \cdot n|}} k_E^2\right)^{1/2}}_{=:T_1}}
       \underbrace{\left( \sum_{E \in \cI} \alpha_E \int_{e_{out}}{|\beta \cdot n_E|}\, 
       \jump{\bar w_h}^2 ds\right)^{1/2}}_{\leq \seminormbeta{w_h}}.
    \end{multline}
    It remains to control $T_1$ using \eqref{eq:kEcontrol} and the definition of $\alpha_E$ given by \eqref{eq:newdefalpha}
    \begin{multline*}
        T_1^2\leq \norm{\beta}_{W^{1,\infty}}^2\sum_{E \in \cI} \frac{\norm{v}^2_{W^{1,1}(E)}}{\alpha_E \int_{e_{out}} |\beta \cdot n|}
        = \norm{\beta}_{W^{1,\infty}}^2 \sum_{E \in \cI} \frac{\tau h}{|E|}\| v \|_{W^{1,1}(E)}^2\\
        \leq \norm{\beta}_{W^{1,\infty}}^2 \sum_{E \in \cI} \tau h \| v \|_{H^{1}(E)}^2
        \leq \norm{\beta}_{W^{1,\infty}}^2 \tau h \| v \|_{H^{1}(\Omega)}^2.
    \end{multline*}
    Here we used that
    \[
    \frac{1}{\abs{E}} \norm{v}_{W^{1,1}(E)}^2 = \frac{1}{\abs{E}}
    \left( \int_E (\abs{v} + \abs{\nabla v})  \right)^2 \le \int_E (\abs{v}^2 + \abs{\nabla v}^2).
    \]
    This concludes the proof.
%
%
\end{proof}  

\subsection{Discrete dissipation}
We begin this section by stating two technical propositions:
\begin{lemma}\label{lemma: (II)}
(i) There holds for every $v_h \in \mathcal{V}^0$ 
\begin{equation}\label{eq: property (II)}
\sum_{e \in \cF} \int_e \omega_e \average{\bar v_h} \beta \cdot \jump{\bar v_h} = 0,
\qquad\text{with\qquad}
\omega_e =
\begin{cases}
       1 & \text{for } e \in \cF_i, 
        \\
        1/2 & \text{for } e \in \cF_b.  
        \end{cases}   
\end{equation}

(ii) There holds for every $v_h, w_h \in \mathcal{V}^0$
\begin{equation}\label{eq:betaid 1} 
     \sum_{e \in \cF} \int_e \beta \cdot \jump{\bar v_h \bar w_h} =0.
 \end{equation}
\end{lemma}

\begin{proof}
(i) We first prove \eqref{eq: property (II)}.
For discrete trial functions $v_h$ we can replace $\beta$-means with actual traces and get
for each inner face $e \in \cF_i$ 
\begin{multline*}
\average{\bar v_h} \beta \cdot \jump{\bar v_h} = 
\average{v_h} \beta \cdot \jump{v_h} = 
\frac 1 2 (v_{E_1} + v_{E_2}) \beta \cdot (v_{E_1} n_{E_1} - v_{E_2} n_{E_1}) \\
= \frac 1 2 (v_{E_1}^2 - v_{E_2}^2) \beta \cdot n_{E_1} 
= \frac 1 2 v_{E_1}^2 \beta \cdot n_{E_1} + \frac 1 2 v_{E_2}^2 \beta \cdot n_{E_2}.
\end{multline*}
For boundary faces $e \in \cF_b$ we get with definition \ref{def: jump and average}
\[
\average{\bar v_h} \beta \cdot \jump{\bar v_h} = 
\average{v_h} \beta \cdot \jump{v_h} = 
v_{E} \: \beta \cdot (v_{E} n_{E})
= v_{E}^2 \: \beta \cdot n_{E}. 
\]
Using the definition of $\omega_e$, we apply the divergence theorem and exploit that $v_h$ is piecewise constant to get
\begin{multline*}
\sum_{e \in \cF} \int_e \omega_e \average{\bar v_h} \beta \cdot \jump{\bar v_h} = 
\frac12 \int_{\partial \Omega}  \beta \cdot n (\bar v_h)^2 +  \sum_{e \in \cF_i} \int_e \left(\tfrac 1 2 v_{E_1}^2 \beta \cdot n_{E_1} + \tfrac 1 2 v_{E_2}^2 \beta \cdot n_{E_2}\right) \\
= \frac 1 2 \sum_{ E \in \cT_h } \int_{\partial E} v_h^2 (\beta \cdot n_E)
= \frac 1 2 \sum_{E \in \cT_h } \int_E \nabla \cdot (\beta v_h^2) = 0.
\end{multline*}

(ii) Now let us proceed to \eqref{eq:betaid 1}. Exploiting that $\nabla \cdot \beta = 0$ and that $\bar v_h$ and $\bar w_h$ are piecewise constant
there holds
\[
0 = \sum_E \int_E (\nabla \cdot \beta) \bar v_h \bar w_h 
= \sum_E \int_{\partial E} (\beta \cdot n) \bar v_h \bar w_h 
= \sum_e \int_e \beta \cdot \jump{\bar v_h \bar w_h}.
\]
This concludes the proof.
\end{proof}

\begin{lemma}\label{lemma: A,B coercivity}
For $A,B \in \mathbb{R}$ and $\alpha \in [0,1]$ there holds
\[
\frac 1 2 A^2 + \alpha B^2 +\alpha AB = \frac 1 2
\left( 
(1-\alpha) A^2 + \alpha B^2 + \alpha (A+B)^2 \right).
\]
\end{lemma}

\begin{proof}
We can rewrite this in matrix form as
\begin{align*}
\frac 1 2 A^2 + \alpha B^2 +\alpha AB &=
\begin{pmatrix}
A \\ B
\end{pmatrix}^T 
\begin{pmatrix}
\frac 1 2 & \frac 1 2 \alpha \\ \frac 1 2 \alpha & \alpha 
\end{pmatrix}
\begin{pmatrix}
A \\ B
\end{pmatrix} \\
&= \frac 1 2 
\left[   
\begin{pmatrix}
A \\ B
\end{pmatrix}^T 
\begin{pmatrix}
(1-\alpha) & 0 \\ 0 & \alpha 
\end{pmatrix}
\begin{pmatrix}
A \\ B
\end{pmatrix}
+ \begin{pmatrix}
A \\ B
\end{pmatrix}^T 
\begin{pmatrix}
\alpha & \alpha \\ \alpha & \alpha 
\end{pmatrix}
\begin{pmatrix}
A \\ B
\end{pmatrix}
\right] \\
&= \frac 1 2
\left( 
(1-\alpha) A^2 + \alpha B^2 + \alpha (A+B)^2 \right) .
\end{align*}
This concludes the proof.
\end{proof}

Now we can show the following result.
\begin{lemma}[Discrete dissipation]\label{lemma: coercivity 2d lin adv P0}
There holds for all $v_h \in \mathcal{V}^0$
\[
\ahDoD(v_h,v_h) = \frac 1 2 \seminormbeta{v_h}^2.
\]
\end{lemma}
\begin{proof}
By definition of $\ahDoD$ and lemma \ref{lemma: (II)}, there holds 
\begin{align*}
    \ahDoD(v_h,v_h) - 0 
    &= \sum_{e \in \cFwoo} \int_e 
    \left( \bar v_h^\uparrow \beta \cdot \jump{\bar v_h} - \omega_e\average{\bar v_h} \beta \cdot \jump{\bar v_h}        \right) \\
    &
    + \sum_{E \in \cI} \int_{e_{out}} \left( \alpha_E \bar v_h^{\uparrow} + (1-\alpha_E) \bar v_{in}   \right) \beta \cdot \jump{\bar v_h} - \omega_e\average{\bar v_h} \beta \cdot \jump{\bar v_h}.
\end{align*}
If $e$ is not a boundary face, i.e., $e 
\in \cF_i \cap \cFwoo$, we can rewrite the first term 
in terms of upwind and downwind traces
as 
\begin{multline}
\int_e 
    \left( \bar v_h^\uparrow \beta \cdot \jump{\bar v_h} - \average{\bar v_h} \beta \cdot \jump{\bar v_h}        \right) =
    \int_e (\bar v_h^\uparrow \beta \cdot (\bar v_h^\uparrow n^\uparrow + \bar v_h^\downarrow n^\downarrow ) - \tfrac12 (\bar v_h^\uparrow +\bar v_h^\downarrow)
    \beta \cdot (\bar v_h^\uparrow n^\uparrow + \bar v_h^\downarrow n^\downarrow )
    \\
    = \int_e\tfrac12 (\bar v_h^\uparrow -\bar v_h^\downarrow) \beta \cdot (\bar v_h^\uparrow n^\uparrow + \bar v_h^\downarrow n^\downarrow )
    = \int_e \tfrac12 (\bar v_h^\uparrow -\bar v_h^\downarrow) \beta \cdot n^\uparrow (\bar v_h^\uparrow - \bar v_h^\downarrow  )
    =  
 \int_e \frac 1 2  \abs{\beta \cdot n} \jump{\bar v_h}^2,
\end{multline}
using that $\beta \cdot n^\uparrow \geq 0$.
If $e$ is a boundary face of a boundary cell $E$, i.e., $e 
\in \cF_b \cap \cFwoo$, we note that we can rewrite the integrand of the first term as (compare also definitions \ref{def: upwind traces} and \ref{def: jump and average})
\begin{equation*}
(\bar v_h^\uparrow - \omega_e\average{\bar v_h})\beta \cdot \jump{\bar v_h} = \left\{
\begin{array}{l@{\:}l@{\:}ll}
(0 - \frac12 \bar v_h) \beta \cdot n \: \bar v_h 
&= \frac12 |\beta \cdot n| \bar v_h^2
&=\frac 1 2  \abs{\beta \cdot n} \jump{\bar v_h}^2
&:\quad e \in \cF_{in}(E),\\[1ex]
(\bar v_h - \frac12 \bar v_h) \beta \cdot n \: \bar v_h 
&= \frac12 |\beta \cdot n| \bar v_h^2
&=\frac 1 2  \abs{\beta \cdot n} \jump{\bar v_h}^2
&:\quad e \in \cF_{out}(E).
\end{array}
\right.
\end{equation*}

For the second term, where $E$ is a stabilized cut cell, we recall that neither $e_{in}$ nor $e_{out}$ is a boundary  face:
\begin{align*}
&\int_{e_{out}} \left( \alpha_E \bar v_{E} + (1-\alpha_E) \bar v_{in}   \right) \beta \cdot \jump{\bar v_h} - \average{\bar v_h} \beta \cdot \jump{\bar v_h}\\
&=\int_{e_{out}} \left( \alpha_E \bar v_{E} - \alpha_E \bar v_{in}   \right) \beta \cdot \jump{\bar v_h} 
+ \bar v_{in} \beta \cdot \jump{\bar v_h}  - \average{\bar v_h} \beta \cdot \jump{\bar v_h}\\
&= \alpha_E \int_{e_{out}} \left(  \bar v_{E} - \bar v_{in}   \right) \beta \cdot \jump{\bar v_h} 
+ \int_{e_{out}}  \bar v_{in} \beta \cdot \jump{\bar v_h}  - \average{\bar v_h} \beta \cdot \jump{\bar v_h}\\
&= \alpha_E \int_{e_{out}} \left(  \bar v_{E} - \bar v_{in}   \right) \beta \cdot n_E (\bar v_{E} - \bar v_{out}) 
+ \int_{e_{out}}  \left( \bar v_{in}  - \average{\bar v_h} \right) \beta \cdot n_E (\bar v_{E} - \bar v_{out})  \\
&= \alpha_E \int_{e_{out}} \left(  \bar v_{E} - \bar v_{in}   \right) \beta \cdot n_E (\bar v_{E} - \bar v_{out}) 
+ \int_{e_{out}}  \left( \bar v_{in}  - \frac{\bar v_{E}+\bar v_{out}}{2} \right) \beta \cdot n_E (\bar v_{E} - \bar v_{out}).
\end{align*}
We will use this convention also in the subsequent computations. This gives us
\begin{align*}
    \ahDoD(v_h,v_h) &= 
    \sum_{e \in \cFwoio} \int_e \frac 1 2  \abs{\beta \cdot n} \jump{\bar v_h}^2 
    +\sum_{E \in \cI} \left( \int_{e_{in}} \frac 1 2 \abs{\beta \cdot n} (\bar v_{in} - \bar v_{E})^2 \right. \\
    &\qquad \left. + \alpha_E \int_{e_{out}} \left(  \bar v_{E} - \bar v_{in}   \right) \beta \cdot n_E (\bar v_{E} - \bar v_{out}) 
+ \int_{e_{out}}  \left( \bar v_{in}  - \frac{\bar v_{E}+\bar v_{out}}{2} \right) \beta \cdot n_E (\bar v_{E} - \bar v_{out}) \right).\\
&= T_1 + \sum_{E \in \cI} \left( T_{\text{2a}} + T_{\text{2b}} + T_{\text{2c}} \right).
    \end{align*}
Comparing this with the definition of the $\beta-$seminorm in \eqref{eq: seminorm beta 2d P0}, we find that the term $T_1$ shows up in the seminorm and that the scaling with the factor $\tfrac 1 2$ is consistent with the claim, i.e., that term is fine. We now take care of the other three terms $T_{\text{2a}}, T_{\text{2b}}, T_{\text{2c}}$. 
We start with $T_{\text{2a}}$ and $T_{\text{2c}}$. 
Using \eqref{eq: betaid combined} on stabilized cut cells, we rewrite integrals over inflow faces as integrals over outflow faces. This gives 
\begin{align*}
T_{\text{2a}} + T_{\text{2c}} &=
\int_{e_{out}} \frac 1 2 \beta \cdot n_E (\bar v_{in} - \bar v_{E})^2
+ \left( \bar v_{in}  - \frac{\bar v_{E}+\bar v_{out}}{2} \right) \beta \cdot n_E (\bar v_{E} - \bar v_{out}) \\
&= \int_{e_{out}} (\beta \cdot n_E) \left[ \frac 1 2 (\bar v_{in} - \bar v_{E})^2
+ \frac 1 2 (\bar v_{in} - \bar v_{E}) (\bar v_{E} - \bar v_{out})
+ \frac 1 2 (\bar v_{in} - \bar v_{out}) (\bar v_{E} - \bar v_{out}) \right] \\
&= \int_{e_{out}} (\beta \cdot n_E) \left[ \frac 1 2 (\bar v_{in} - \bar v_{E})(\bar v_{in} - \bar v_{E} + \bar v_{E} - \bar v_{out})
+ \frac 1 2 (\bar v_{in} - \bar v_{out}) (\bar v_{E} - \bar v_{out}) \right] \\
&= \int_{e_{out}} (\beta \cdot n_E) \left[ \frac 1 2 (\bar v_{in} - \bar v_{E})(\bar v_{in}  - \bar v_{out})
+ \frac 1 2 (\bar v_{in} - \bar v_{out}) (\bar v_{E} - \bar v_{out}) \right]  \\
&= \int_{e_{out}} (\beta \cdot n_E) \frac 1 2 (\bar v_{in} - \bar v_{out})(\bar v_{in} - \bar v_{E} + \bar v_{E} - \bar v_{out}) \\
&= \int_{e_{out}} (\beta \cdot n_E) \frac 1 2 (\bar v_{in} - \bar v_{out})^2.
\end{align*}
Now we still need to consider $T_{\text{2b}}$ given by
\[
\alpha_E \int_{e_{out}} \left(  \bar v_{E} - \bar v_{in}   \right) \beta \cdot n_E (\bar v_{E} - \bar v_{out}) .
\]
There holds
\begin{align*}
(  \bar v_{E} - \bar v_{in}) (\bar v_{E} - \bar v_{out}) 
&= (\bar v_{E} - \bar v_{in}) (\bar v_{E} - \bar v_{in} + \bar v_{in} - \bar v_{out}) \\
&= (\bar v_{E} - \bar v_{in})^2 + \underbrace{(\bar v_{E} - \bar v_{in})}_{=:B} \underbrace{(\bar v_{in} - \bar v_{out})}_{=:A} \\
&= B^2 + AB.
\end{align*}
Note that $A+B = \bar v_{E} - \bar v_{out}$.
Therefore
\[
T_{\text{2a}} + T_{\text{2c}} + T_{\text{2b}} =
\int_{e_{out}} (\beta \cdot n_E) \left(\frac 1 2 A^2 + \alpha_E B^2 +\alpha_E AB\right).
\]
Applying lemma \ref{lemma: A,B coercivity} gives
\begin{align*}
T_{\text{2a}} + T_{\text{2c}} + T_{\text{2b}} &= 
\frac 1 2 \int_{e_{out}} (\beta \cdot n_E) ((1-\alpha_E) A^2 + \alpha_E B^2 + \alpha_E (A+B)^2 ) \\
&= \frac 1 2 \int_{e_{out}} (\beta \cdot n_E) \left[ 
(1 - \alpha_E)(\bar v_{in}-{\bar v_{out}})^2 + \alpha_E (\bar v_{E}-\bar v_{in})^2 + \alpha_E(\bar v_{E}- \bar v_{out})^2
\right].
\end{align*}
Using now equation \eqref{eq: betaid combined} again to rewrite the second term as integral over face $e_{in}$, these terms match exactly, up to the factor $\tfrac 1 2$, the terms in the second line in the definition of the seminorm. This concludes the proof. 
\end{proof}

\subsection{Boundedness}



\begin{lemma}[Boundedness I]\label{lem:(3.33)}
    For all $v \in \VhStern$ and $w_h \in \mathcal{V}^0$ it holds that
    \[
   \abs{ \ahDoD (v,w_h)}  \leq \normuwbStern{v} |w_h|_\beta.
    \]
%
%
%
%
%
\end{lemma}

\begin{proof}
  We have from \eqref{eq: def DoD upw} by definition of $\ahDoD$
  \begin{align*} 
\ahDoD (v,w_h)
&= \sum_{e \in \cFwoo} \int_e 
     \bar v^\uparrow \beta \cdot \jump{\bar w_h} + \sum_{E \in \cI} \int_{e_{out}(E)} \left( \alpha_E \bar v_{E} + (1-\alpha_E) \bar v_{in}   \right) \beta \cdot \jump{\bar w_h} \\
&= \underbrace{\sum_{e \in \cFwoio} \int_e  \bar v^\uparrow \beta \cdot \jump{\bar w_h}}_{=:T_1} \\
& \quad + \underbrace{\sum_{E \in \cI}\left[\int_{e_{in}} (\beta \cdot n_E) (\bar w_E - \bar w_{in} ) \bar v_{in} 
+ \int_{e_{out}} (\beta \cdot n_E) (\bar w_E - \bar w_{out} )( (1-\alpha_E)\bar v_{in} + \alpha_E \bar v_E)\right]}_{=:T_2}.\\
  \end{align*}

In the last step, we extracted the boundary integrals over inflow faces of stabilized cut cells from the first sum and added it to the second, and explicitly wrote out what $\jump{\bar w_h}$ means on inflow and outflow faces of stabilized cells, respectively.

We now rewrite parts of $T_2$ by using \eqref{eq: betaid combined} to get
  \begin{multline*}
  \int_{e_{in}} (\beta \cdot n_E) (\bar w_E - \bar w_{in} ) \bar v_{in} =
  \int_{e_{in}} (\beta \cdot n_E) (\bar w_E - \bar w_{in} ) \alpha_E \bar v_{in} + 
  \int_{e_{in}} (\beta \cdot n_E) (\bar w_E - \bar w_{in} ) (1-\alpha_E)\bar v_{in} \\
  =\int_{e_{in}} (\beta \cdot n_E) (\bar w_E - \bar w_{in} ) \alpha_E \bar v_{in} -
  \int_{e_{out}} (\beta \cdot n_{E}) (\bar w_E - \bar w_{in} ) (1 - \alpha_E) \bar v_{in}
  \end{multline*}
  for each stabilized cut cell $E$, which leads to
\begin{align*}
  T_2 &= \sum_{E \in \cI}\left[
  \int_{e_{in}} (\beta \cdot n_E) (\bar w_E - \bar w_{in} ) \alpha_E \bar v_{in} -
  \int_{e_{out}} (\beta \cdot n_{E}) (\bar w_E - \bar w_{in} ) (1 - \alpha_E) \bar v_{in} \right.\\
&\qquad + \left. \int_{e_{out}} (\beta \cdot n_E) (\bar w_E - \bar w_{out} )(1-\alpha_E)\bar v_{in} + \int_{e_{out}} (\beta \cdot n_E) (\bar w_E - \bar w_{out} )\alpha_E \bar v_E   \right] \\
&= \sum_{E \in \cI}\left[\int_{e_{in}} (\beta \cdot n_E) (\bar w_E - \bar w_{in} ) \alpha_E \bar v_{in} \right.\\
&\qquad \left. + \int_{e_{out}} (\beta \cdot n_E) (\bar w_{in} - \bar w_{out} ) (1-\alpha_E)\bar v_{in} + \int_{e_{out}} (\beta \cdot n_E) (\bar w_E - \bar w_{out}) \alpha_E \bar v_E)\right].
  \end{align*}
Summing up $T_1$ and $T_2$ and writing out jumps in $w_h$ over non-stabilized faces as well (denoting by $E_1$ and $E_2$ the 2 cells that share the face $e$) we get
  \begin{align*}
  T_1 + T_2 
&= \sum_{e \in \cFwoio} \int_e  \bar v^\uparrow \beta \cdot n_{E_1} (\bar w_{E_1} - \bar w_{E_2}) + \sum_{E \in \cI}\left[\int_{e_{in}} (\beta \cdot n_E) (\bar w_E - \bar w_{in} ) \alpha_E \bar v_{in} \right.\\
&\qquad \left. + \int_{e_{out}} (\beta \cdot n_E) (\bar w_{in} - \bar w_{out} ) (1-\alpha_E)\bar v_{in} + \int_{e_{out}} (\beta \cdot n_E) (\bar w_E - \bar w_{out}) \alpha_E \bar v_E)\right].
  \end{align*}
By using Cauchy-Schwarz this gives the bound
\[
T_1 + T_2 \le \sqrt{T_3} \seminormbeta{w_h}
\]
with
\[
T_3 = \sum_{e \in \cFwoio} \int_e \abs{\beta \cdot n} (\bar v^\uparrow)^2 + \sum_{E \in \cI} \left[ \alpha_E \int_{e_{in}} |\beta \cdot n_E| \bar v_{in}^2 + \alpha_E \int_{e_{out}} |\beta \cdot n_E| \bar v_E^2 + (1-\alpha_E) \int_{e_{out}} |\beta \cdot n_E| \bar v_{in}^2 \right].
\]
It remains to bound $T_3$. Using \eqref{eq: betaid combined} again and the fact that $\bar v_{in} = \vweighted$ is just a constant 
we have that
\[
(1 - \alpha_E) \int_{e_{out}} |\beta \cdot n_{E}| \bar v_{in}^2 = (1 - \alpha_E)
\int_{e_{in}} |\beta \cdot n_E|  \bar v_{in}^2 .
\]
Summarizing the terms involving $\bar v_{in}^2$ using $(1-\alpha_E)+\alpha_E=1$ gives us the bound
\[
    T_3 = \sum_{e \in \cFwoio} \int_e \abs{\beta \cdot n} (\bar v^\uparrow)^2 + \sum_{E \in \cI} \left[ \int_{e_{in}} |\beta \cdot n_E| \bar v_{in}^2 + \alpha_E \int_{e_{out}} |\beta \cdot n_E| \bar v_E^2  \right] \le \normuwbStern{v}^2,
\]
which finishes the proof.

 \end{proof}

  \begin{lemma}[Boundedness II]\label{lem:1bound}
Given any $v_h, w_h \in \mathcal{V}^0$ the following two estimates hold
\begin{align}\label{eq:lem:1bound:1}
&\text{\normalfont(i)}&
\ahDoD(v_h, w_h)
&\le \sqrt{\frac{\Ctr}{h}} |v_h|_\beta \|w_h \|_{L^2(\Omega)},
\text{ and}\\
\label{eq:lem:1bound:2}
&\text{\normalfont(ii)}&
\| \AhDoD v_h \|_{L^2(\Omega)} 
&\leq \sqrt{\frac{\Ctr}{h}} \seminormbeta{v_h},
\end{align}
with $\Ctr$ from lemma \refwithname{lem: inv est const}.
\end{lemma}

\begin{proof}
Statement (ii), i.e., estimate \eqref{eq:lem:1bound:2}, follows with \eqref{eq: operator AhDoD}
directly from \eqref{eq:lem:1bound:1}. So we proceed with proving statement (i). To do so, 
we are going to use equation \eqref{eq:betaid 1}, the 
algebraic identity
 \begin{equation}\label{eq: jump of prod}
     \jump{\bar v_h \bar w_h} = \bar v_h^\uparrow \jump{\bar w_h} + {\bar w_h^{\downarrow}}\jump{\bar v_h},
 \end{equation}
 and \eqref{eq: betaid combined} for $E$ being a stabilized cut cell to rewrite an integral over the inflow face as an integral over the outflow face and vice versa.
Then, 
\begin{align*}
    \ahDoD&(\bar v_h, \bar w_h) \overset{\mathclap{\eqref{eq:betaid 1}}}{=} \ahDoD(\bar v_h, \bar w_h) - \sum_{e \in \cF} \int_e \beta \cdot \jump{\bar v_h \bar w_h}\\
    &\overset{\mathclap{\eqref{eq: def DoD upw}}}{=} 
    \sum_{e \in \cFwoo} 
    \int_e \left( \bar v_h^\uparrow \beta \cdot \jump{\bar w_h} - \beta \cdot \jump{\bar v_h \bar w_h} \right) + 
    \sum_{E \in \cI}\int_{e_{out}} \left( \alpha_E \bar v_{E} + (1-\alpha_E) \bar v_{in}   \right) \beta \cdot \jump{\bar w_h} 
    - \beta \cdot \jump{\bar v_h \bar w_h} \\
    &\overset{\mathclap{\eqref{eq: jump of prod}}}{=} 
    - \sum_{e \in \cFwoo} \int_e  {\bar w_h^{\downarrow}}\beta \cdot \jump{\bar v_h}\\
    &\quad +  \sum_{E \in \cI} \int_{e_{out}} (\alpha_E \bar v_E + (1- \alpha_E)\bar v_{in} )(\beta \cdot n_E) (\bar w_E - \bar w_{out})
       - (\beta \cdot n_E) (\bar v_E \bar w_E - \bar v_{out} \bar w_{out} )\\
&= -\sum_{e \in \cFwoio} \int_e  {\bar w_h^{\downarrow}}\beta \cdot \jump{\bar v_h}  
- \sum_{E \in \cI} \int_{e_{in}} \bar w_E (\beta \cdot n_E) (\bar v_E - \bar v_{in})\\
&\phantom{=}~ + \sum_{E \in \cI} \int_{e_{out}} \bar{w}_E (\beta \cdot n_E) [ (1- \alpha_E) \bar{v}_{in} - (1- \alpha_E)\bar{v}_E] - \bar{w}_{out} (\beta \cdot n_E) [ \alpha_E \bar{v}_E + (1-\alpha_E) \bar{v}_{in} - \bar{v}_{out}]\\
&\overset{\mathclap{\eqref{eq: betaid combined}}}{=} -\sum_{e \in \cFwoio} \int_e  {\bar w_h^{\downarrow}}\beta \cdot \jump{\bar v_h}  
+ \sum_{E \in \cI} \int_{e_{in}} \bar w_E (\beta \cdot n_E) (\bar v_{in} - \bar v_E)
- \sum_{E \in \cI} (1- \alpha_E) \int_{e_{in}} \bar{w}_E (\beta \cdot n_E) [ \bar{v}_{in} - \bar{v}_E]
\\
&\phantom{=}~ - \sum_{E \in \cI} \int_{e_{out}} \bar{w}_{out} (\beta \cdot n_E) [ \alpha_E \bar{v}_E 
- \alpha_E \bar v_{out} + (1-\alpha_E) \bar{v}_{in} - (1-\alpha_E) \bar{v}_{out}]\\
&= - \sum_{e \in \cFwoio} \int_e  {\bar w_h^{\downarrow}}\beta \cdot \jump{\bar v_h} 
+  \sum_{E \in \cI} \int_{e_{out}} (\beta \cdot n_E) (1- \alpha_E) (\bar v_{out} - \bar v_{in}) \bar w_{out}\\
&\quad - \sum_{E \in \cI} \alpha_E \biggl(\int_{e_{out}}  \bar w_{out} \beta \cdot \jump{\bar v_h}  + \int_{e_{in}} \bar w_E \beta \cdot \jump{\bar v_h} \biggr).
\end{align*}

By definition \ref{def: seminorm beta 2d P0} of $\seminormbeta{\cdot}$ and the Cauchy-Schwarz inequality we can bound this by
\begin{align*}
    \ahDoD(\bar v_h, \bar w_h) & \le |\bar v_h|_\beta \biggl(\sum_{e \in \cFwoio} \int_e |\beta \cdot n | {\bar w_h^{\downarrow^2}}
+  \sum_{E \in \cI} \int_{e_{out}} |\beta \cdot n_E| (1- \alpha_E) \bar w_{out}^2  \\ & \quad + \sum_{E \in \cI}\alpha_E \biggl [ \int_{e_{out}}  |\beta \cdot n| \bar w_{out}^2 + \int_{e_{in}} |\beta \cdot n| \bar w_E^2 \biggr] \biggr)^{\frac{1}{2}}\\
&\le |\bar v_h|_\beta \biggl(\sum_{e \in \cFwoio} \int_e |\beta \cdot n | {\bar w_h^{\downarrow^2}}
+  \sum_{E \in \cI} \int_{e_{out}} |\beta \cdot n_E| \bar w_{out}^2 + \sum_{E \in \cI} \alpha_E  \int_{e_{in}} |\beta \cdot n| \bar w_E^2 \biggr)^{\frac{1}{2}}.
\end{align*}
Let us first consider the first two sums. These contain only functions with support on unstabilized cells $E \not \in \cI$.
Sorting the contributions that belong to a cell $E \not \in \cI$, we observe that for the first sum this corresponds to the sum over inflow faces of cell $E$, i.e., faces with $\beta \cdot n_E < 0$. 
Including the second sum, for any $E \not \in \cI$ that is an outflow neighbor of a stabilized cut cell there holds $\bar w_{out}=\bar w_h = \bar w_h^{\downarrow}$, i.e., viewed from cell $E$'s perspective, it is just a normal inflow face and can be combined with the inflow terms from the first sum. 
Using \eqref{eq:inv est const unstab} from lemma \ref{lem: inv est const}, we establish the bound
\[
\sum_{e \in \cFwoio} \int_e |\beta \cdot n | {\bar w_h^{\downarrow^2}}
+  \sum_{E \in \cI} \int_{e_{out}} |\beta \cdot n_E| \bar w_{out}^2 \le \frac{ \Ctr[1]}{h} \sum_{E \not \in \cI} \| \bar w_h \|_{L^2(E)}^2.
\]
The third sum only contains functions with support on stabilized cut cells. Again the contributions stem only from inflow faces. Thus using \eqref{eq:inv est const stab} we can bound it by
\[
\sum_{E \in \cI} \alpha_E  \int_{e_{in}} |\beta \cdot n| \bar w_E^2 \le \frac{\Ctr[2]}{h} \sum_{E \in \cI} \| \bar w_h \|_{L^2(E)}^2.
\]
Putting everything together we arrive at
\[
\ahDoD(\bar v_h, \bar w_h) \le \frac{\sqrt{\max(\Ctr[1], \Ctr[2])}}{\sqrt{h}} |\bar v_h|_\beta \|\bar w_h \|_{L^2(\Omega)}
= \sqrt{\frac{\Ctr}{h}} |\bar v_h|_\beta \|\bar w_h \|_{L^2(\Omega)}
\]
which is the desired result.
\end{proof}

\section{Proof of main result}\label{sec: proof of main result}
To derive our final error estimate we need the following discrete Gronwall lemma:
\begin{lemma}\label{gronwall}
    Let non-negative sequences $\{ c_n\}_{n \ge 0} ,\{ w_n\}_{n \ge 0} \subset \mathbb{R}$ be given such that  for constants $K, \dt >0$ there holds $0 \le w_{n+1} \le (1+\dt K)w_n + \dt c_n$. 
Then, for all $N \in \mathbb{N}$ there holds 
\[
w_N \le w^0 e^{K \dt N} + \max_{0\le n \le N-1} c_n (e^{\dt K N}-1)/K.
\]
\end{lemma}

\begin{theorem}[Discrete stability]\label{thm:discrete stability}
Let $u \in C^2([0,T], L^2(\Omega)) \cap C^0([0,T], H^1(\Omega))$ be the exact solution of \eqref{eq: lin adv 2d}. Let $u_h^n$ for $n=1,\ldots,N$ be the solution of \eqref{eq: scheme 2d with stab fully discrete} and let the assumptions on the triangulation and the velocity field $\beta$ be as described in section \ref{sec: prelim}.
Fix any $0 <\varepsilon < \tfrac 1 2.$
Let the time step satisfy the CFL-condition 
 \begin{equation}\label{eq: CFL cond}
 \dt \leq \kappa h, \quad \kappa = \frac{1-2\varepsilon}{(1+\varepsilon)\max(4\norm{\beta}_{\infty},\tau^{-1})},
 \end{equation}
 i.e., with $\kappa$ independent of $h$ and of the size of cut cells. Let $N\dt = T$. Define 
\[
\xihn \coloneqq u_h^n - \Pi_h u^n \: \text{(discrete error)}, \quad 
\xipin \coloneqq u^n - \Pi_h u^n \: \text{(projection error)}
\]
with $\Pi_h: \VhStern \rightarrow \mathcal{V}^0$ denoting the $L^2$-projection.
Then there are constants $\genC[i] > 0$ such that 
\begin{equation*}
\norm{\xihnplus}^2 + \varepsilon\dt \seminormbeta{\xihn}^2
\le \left(1 + \frac \dt T\right)\norm{\xihn}^2  
+ \dt (1+\epsinv)\left[  
 \genC[1]
\normuwbStern{\xipin} ^2 
+ \left(T+\genC[3]\dt \right) \norm{\theta^n}^2
+ \genC[4] h
\|u^n\|^2_{H^{1}(\Omega)}\right]
\end{equation*}
with 
\begin{equation*}
    \theta^n = \frac{1}{\dt} \int_{t^n}^{t^{n+1}} (t^{n+1}-t)d_t^2 u(t) \: dt.
\end{equation*}
The generic constants $\genC[i]$ are independent of $h,\Delta t$, the size of small cut cells, the exact solution, the final time $T$, and the choice of $\varepsilon$. 
\end{theorem}

\begin{remark}
For the case $\norm{\beta}_{\infty}=1$ and the choices $\tau=1$ and $\varepsilon = \tfrac{1}{14}$,
\eqref{eq: CFL cond} corresponds to the CFL condition $\kappa = \frac{1}{5}$. In the limit for $\varepsilon \to 0$, it even goes down to $\kappa \approx \frac 1 4$. This is not far away from the CFL condition that one would use in practice, given that the larger cut cells are not stabilized. 
\end{remark}

\begin{remark}
For better readability, dependencies on the velocity field $\beta$ and the parameter $\tau$ are included in the constant $C$. This holds for the results in both theorem \ref{thm:discrete stability} and in theorem \ref{thm:main}. We note though that in theorem \ref{thm:main} the dependence of the constant on $\beta$ is at most quadratic and on $\frac{1}{\tau}$ it is at most of order four.
\end{remark}

\begin{remark}\label{remark: Taylor}
The term $\theta^n$ corresponds to the remainder term of a Taylor series expansion of the exact solution of the form
$u^{n+1} = u^n + \Delta t d_t u^n + \Delta t \theta^n$.
\end{remark}

\begin{proof}
Similarly to $\AhDoD$, we introduce operators $\Ahupw$ and $J_h$ where
$\AhDoD = \Ahupw + J_h$. (We misuse notation here and use $J_h$ to denote the operator $\VhStern \rightarrow \mathcal{V}^0$ induced by the bilinear form $J_h(\cdot,\cdot)$ as well as the bilinear form itself.)
Then $\Ahupw = \AhDoD - J_h$. We also introduce the operator $L_h: (0, T) \to \mathcal{V}^0$ induced by the right hand side $l_h(\cdot,\cdot)$, compare \eqref{eq: ah lin adv deg p l_h}. 

Doing Taylor expansion as indicated in remark \ref{remark: Taylor}, applying the $L^2$ projection, and using the consistency of the upwind DG bilinear form $a_h^{upw}(\cdot,\cdot)$ to replace $\Pi_h d_t u^n$ by $-\Ahupw u^n - L_h^n$, we get
\begin{align}\label{eq: Aug 26 eq 1}
\begin{split}
    \Pi_h u^{n+1} &= \Pi_h u^n - \dt \Ahupw u^n - \dt L_h^n + \dt \Pi_h \theta^n \\
    &= \Pi_h u^n - \dt \AhDoD u^n - \dt L_h^n + \dt \Pi_h \theta^n +\dt J_h(u^n).
    \end{split}
\end{align}
%
Further, from the discrete scheme, we get
\begin{equation}\label{eq: Aug 26 eq 2}
u_h^{n+1} = u_h^n - \dt \AhDoD u_h^n - \dt L_h^n.
\end{equation}
Subtracting \eqref{eq: Aug 26 eq 1} from \eqref{eq: Aug 26 eq 2} results in
\[
\xihnplus = \xihn - \dt \AhDoD \underbrace{(u_h^n - u^n)}_{\xihn - \xipin} - \dt \Pi_h \theta^n - \dt J_h(u^n).
\]
Note that
\[
\xipin - \xihn = (u^n - \Pi_h u^n) - (u_h^n - \Pi_h u^n) = u^n - u_h^n.
\]
This gives the error equation
\begin{equation}\label{eq: error eq 2d}
\xihnplus = \xihn - \dt \AhDoD \xihn + \dt \AhDoD\xipin - \dt \Pi_h\theta^n - \dt J_h(u^n).
\end{equation}
For any $a, b \in \mathbb{R}$, there holds
\[
\frac 1 2 b^2 - \frac 1 2 a^2 - \frac 1 2 (b-a)^2 = 
\frac 1 2 b^2 - \frac 1 2 a^2 - \left(\frac 1 2 b^2 - ab + \frac 1 2 a^2\right) = - a^2 + ab.
\]
With $a=\xihn$ and $b = \xihnplus$, this implies
\[
(\xihnplus, \xihn)_{L^2}- \norm{\xihn}^2 = 
\frac 1 2 \norm{\xihnplus}^2 - \frac 1 2 \norm{\xihn}^2 - \frac 1 2 \norm{\xihnplus - \xihn}^2,
\]
where we write $\| \cdot\|$ for the $L^2(\Omega)$ norm for brevity.

Testing the error equation \eqref{eq: error eq 2d} with $\xihn$ and integrating over the domain, and using the above manipulation results in 
\begin{align*}
\frac 1 2 \norm{\xihnplus}^2 - \frac 1 2 \norm{\xihn}^2 - \frac 1 2 \norm{\xihnplus - \xihn}^2 = 
- \dt \ahDoD(\xihn,\xihn) + \dt \ahDoD(\xipin,\xihn) - \dt (\theta^n,\xihn) - \dt J_h(u^n,\xihn).
\end{align*}
Using now results from above we can estimate the terms on the right hand side as follows:
\begin{itemize}
\item from lemma \refwithname{lemma: coercivity 2d lin adv P0}, we get
$\ahDoD(\xihn,\xihn) = \frac 1 2 \seminormbeta{\xihn}^2 $;
\item from lemma \refwithname{lem:(3.33)}, we have
$\abs{\ahDoD(\xipin,\xihn)}\le \normuwbStern{\xipin} \seminormbeta{\xihn}$;
\item from lemma \refwithname{lemma: consistency_new}, we get
$\abs{J_h(u^n,\xihn)} \le \sqrt{\tau  h} \norm{\beta}_{W^{1,\infty}} \|u^n\|_{H^{1}(\Omega)} \seminormbeta{\xihn}$; in the following we will replace $\sqrt{\tau}\norm{\beta}_{W^{1,\infty}}$ by a generic constant $C$ for brevity.
\end{itemize}
In summary, this gives
\begin{multline*}
\frac 1 2 \norm{\xihnplus}^2 \le \frac 1 2 \norm{\xihn}^2 + \frac 1 2 \norm{\xihnplus - \xihn}^2
- \frac{\dt}{2} \seminormbeta{\xihn}^2 \\
+ \dt \normuwbStern{\xipin} \seminormbeta{\xihn} 
+ \dt \norm{\theta^n} \norm{\xihn}
+ \dt C\sqrt{h}  \|u^n\|_{H^{1}(\Omega)}\seminormbeta{\xihn}.
\end{multline*}
We now apply Young's inequality $ab \le \frac{a^2}{2\delta} + \frac{\delta b^2}{2}$ with $\delta=\frac{\varepsilon}{2}$ and $\delta=\frac{1}{T}$ 
respectively to get
\begin{multline*}
\frac 1 2 \norm{\xihnplus}^2 \le \frac 1 2 \norm{\xihn}^2 + \frac 1 2 \norm{\xihnplus - \xihn}^2 
- \frac{\dt}{2} \seminormbeta{\xihn}^2 \\
+ \dt \left(\frac{1}{\varepsilon} \normuwbStern{\xipin} ^2 + \frac{\varepsilon}{4} \seminormbeta{\xihn}^2 \right)
+ \dt \left( \frac{T}{2} \norm{\theta^n}^2 + \frac{1}{2 T} \norm{\xihn}^2 \right)
+ \dt \left( \frac{1}{\varepsilon} C \: h \|u^n\|^2_{H^{1}(\Omega)}  + \frac{\varepsilon}{4} \seminormbeta{\xihn}^2\right).
\end{multline*}
We collect 
all terms on the right hand side that contain $\seminormbeta{\xihn}^2$ and bring them to the left. This gives
\begin{multline}\label{eq: help1tria}
\frac 1 2 \norm{\xihnplus}^2 + (1-\varepsilon)\frac{\dt}{2}  \seminormbeta{\xihn}^2
\le \left(\frac{1}{2} + \frac{\dt }{2 T}\right)\norm{\xihn}^2 
+ \frac 1 2 \norm{\xihnplus - \xihn}^2 \\
+ \frac{1}{\varepsilon} \dt  \normuwbStern{\xipin} ^2 
+ \frac{\dt T}{2} \norm{\theta^n}^2
+ \frac{1}{\varepsilon} C\dt \: h  \|u^n\|^2_{H^{1}(\Omega)} .
\end{multline}
In the following we will derive an estimate for $\norm{\xihnplus - \xihn}^2 $. From the error equation 
\eqref{eq: error eq 2d}, we get
\[  
\xihnplus - \xihn = - \dt \AhDoD \xihn + \dt \AhDoD\xipin - \dt \Pi_h \theta^n - \dt J_h(u^n).
\]
Taking the $L^2$ norm on both sides and using
 the fact that (which follows essentially from Young's inequality with $\delta = \varepsilon^{-1}$)
\[ \norm{f_1+f_2}^2 \le (1+ \varepsilon) \norm{f_1}^2 + (1+ \varepsilon^{-1})\norm{f_2}^2 
\qquad \forall \varepsilon>0,\ \forall f_i \in L^2(\Omega)\]
so that
\[
\norm{f_1+f_2+f_3+f_4}^2 \le ( 1+ \varepsilon)\norm{f_1}^2 +3(1+ \varepsilon^{-1}) (\norm{f_2}^2 +\norm{f_3}^2 + \norm{f_4}^2),
\quad  \forall \varepsilon>0,\ \forall f_i \in L^2(\Omega),
\]
we get
\begin{equation*}  
\norm{\xihnplus - \xihn}^2
\le (1+\varepsilon) \dt^2 \underbrace{\norm{\AhDoD \xihn}^2}_{=:T_1} 
+   3 \underbrace{(1+ \varepsilon^{-1})}_{=:\taueps} \dt^2 \big(\underbrace{\norm{\AhDoD\xipin}^2}_{=:T_2}
+   \norm{ \theta^n}^2
+   \underbrace{\norm{J_h(u^n)}^2}_{=:T_3}\big),
\end{equation*}
where we also used $\norm{\Pi_h \theta^n}^2 \le \norm{\theta^n}^2$ from the stability of the $L^2$ projection.
In the following we will abbreviate $\taueps := 1+\frac 1 \varepsilon$.
We now individually bound the terms on the right hand side:
\begin{itemize}
\item $T_1$: by lemma \refwithname{lem:1bound} there holds
\[  \norm{\AhDoD \xihn}^2 \leq \frac{\Ctr}{h} \seminormbeta{\xihn}^2, \]
\item $T_2$: by lemma \refwithname{lem:(3.33)} and the inverse estimate from lemma \ref{lemma: inverse eq betaseminorm}, there holds
\[  
(\AhDoD v,w_h) \le  \normuwbStern{v} \seminormbeta{w_h}
\le 2\sqrt{\frac{\Ctr}{h}} \normuwbStern{v} \norm{w_h}
\qquad \forall v \in \VhStern, w_h \in \mathcal{V}^0%
.
\]
This then implies
\[
\norm{\AhDoD v} \le \frac{C}{\sqrt{h}} \normuwbStern{v}
\qquad \forall v \in \VhStern,
\]
where we replaced the constant $\Ctr$ by the generic constant $C$ for brevity,
and therefore
\[
\norm{\AhDoD \xipin}^2 \le \frac{C}{h} \normuwbStern{\xipin}^2;
\]
\item $T_3$: 
applying lemma \refwithname{lemma: consistency_new} and lemma
\refwithname{lemma: inverse eq betaseminorm}
we get
\[
(J_h(u^n),w_h) \le C \sqrt{h} \|u^n\|_{H^{1}(\Omega)} \seminormbeta{w_h}
\le C  \|u^n\|_{H^{1}(\Omega)}\norm{w_h},
\]
which implies
\[
\norm{J_h(u^n)}^2 \le C \|u^n\|^2_{H^{1}(\Omega)}.
\]
\end{itemize}
In total, this gives
\begin{equation}\label{eq: help2tria}  
\norm{\xihnplus - \xihn}^2 \le  (1+\varepsilon) \frac{\Ctr}{h}
\dt^2 \seminormbeta{\xihn}^2 
+ \taueps \dt^2 \left(
     \frac{C}{h} \normuwbStern{\xipin}^2 \\
+    C\norm{\theta^n}^2 
+    C\|u^n\|^2_{H^{1}(\Omega)}\right).
\end{equation}
Using now \eqref{eq: help2tria} in \eqref{eq: help1tria} gives
\begin{multline}
\frac 1 2 \norm{\xihnplus}^2 + \frac{(1-\varepsilon)\dt}{2} \seminormbeta{\xihn}^2
\le \left(\frac{1}{2} + \frac{\dt }{2 T}\right) \norm{\xihn}^2 + \frac{1+\varepsilon}{2}\frac{\Ctr}{h}
\dt^2 \seminormbeta{\xihn}^2 \\
+ \dt \left[  
\left( \varepsilon^{-1}+\taueps C\frac{\dt}{h}\right) \normuwbStern{\xipin} ^2 
+ \left(\frac{T}{2}+\taueps C\dt\right) \norm{\theta^n}^2
+ \taueps C (h+\dt)  \|u^n\|^2_{H^{1}(\Omega)} .\right]
\end{multline}
Choosing a sufficiently small time step length $\dt$ we can guarantee
\[
\frac{1+\varepsilon}{2} \frac{\Ctr}{h} \dt^2 \seminormbeta{\xihn}^2 \overset{!}{\le} \left(\frac{1-2\varepsilon}{2}\right)\dt \seminormbeta{\xihn}^2,
\]
which implies the condition 
\begin{equation}\label{eq: CFL cond proof}
\dt \overset{!}{\le} \frac{1-2\varepsilon}{1+\varepsilon} \frac{h}{\Ctr}.
\end{equation}
As all involved constants are independent of $h$, this gives a time step restriction of the form $\Delta t \le \kappa h$ with $\kappa = \frac{1-2\varepsilon}{(1+\varepsilon)\Ctr} = \frac{1-2\varepsilon}{(1+\varepsilon)\max(4\norm{\beta}_{\infty},\tau^{-1})}$, as assumed in 
\eqref{eq: CFL cond}.
Bounding $\frac{\Delta t}h\le\kappa$ and, equivalently, $\Delta t\le\kappa h$ yields
\begin{multline*}
\frac 1 2 \norm{\xihnplus}^2 + \varepsilon\frac{\dt}{2} \seminormbeta{\xihn}^2
\le \left(\frac{1}{2} + \frac{\dt}{2 T}\right) \norm{\xihn}^2
+ \frac \dt 2  \left[ \taueps
C\normuwbStern{\xipin} ^2 
+ (T+ \taueps C\dt ) \norm{\theta^n}^2
+ \taueps C h\|u^n\|^2_{H^{1}(\Omega)} \right].
\end{multline*}
Multiplying the whole equation by $2$ and factoring out $\taueps$ gives the desired result.
\end{proof}

\begin{theorem}[Fully discrete a-priori error estimate]\label{thm:main}
 Let all assumptions from theorem \ref{thm:discrete stability} hold true. Let $\dt \le \kappa h$ and let $h$ be sufficiently small (so that higher-order terms can be neglected).
 Then, there exists a constant $C>0$ (independent of $h,\dt$, the size of small cut cells, the final time $T$, and the solution $u$) 
 such that for all $M$ with $M \dt \leq T$ there holds
\begin{multline}
 \label{eq:fully-discrete-error-estimate}
    \norm{u(t^M,\cdot)-u_h^M}^2_{L^2(\Omega)} + \sum_{n=0}^{M-1} \varepsilon\dt \seminormbeta{u(t^n,\cdot) - u_h^n}^2 \\
    \le (C + T e^{\frac{\dt M}{T} }C) (1+\epsinv) h
    \norm{u}^2_{C^0([0,T],H^1(\Omega))} + (1+\epsinv) T^2(e^{\frac{\dt M}{T} }-1)\dt^2 \norm{u}^2_{C^2([0,T],L^2(\Omega))}.
\end{multline}
\end{theorem}

\begin{proof}
We start with the result from theorem \ref{thm:discrete stability} and apply lemma \ref{gronwall}
with $w_n := \norm{\xihn}^2 + \sum_{j=0}^{n-1} \varepsilon\dt \seminormbeta{\xi_h^j}^2$. As
\[
w_{n+1} - w_n = \norm{\xihnplus}^2 - \norm{\xihn}^2+ \varepsilon\dt \seminormbeta{\xihn}^2
\]
the assumptions for lemma \ref{gronwall} (of $w_{n+1}-w_n \le \dt K w_n + \dt c_n$) are satisfied with $K=\frac{1}{T}$ and
\[  
c_n =  (1+\epsinv)\left[C \normuwbStern{\xipin} ^2 
+ (T+C\dt) \norm{\theta^n}^2
+ C h \|u^n\|^2_{H^{1}(\Omega)} \right]
\]
due to theorem \ref{thm:discrete stability}.
We therefore get
\begin{equation*}
    \norm{\xi_h^M}^2 + \sum_{j=0}^{M-1} \varepsilon\dt
    \seminormbeta{\xi_h^j}^2
    \le \norm{\xi_h^0}^2 e^{\frac{\dt M}{T}} + T \max_{0\le n \le M-1} c_n (e^{\frac{\dt M}{T} }-1).
\end{equation*}
By definition, $\xi_h^0 = u_h^0 - \Pi_h u_0 = 0$.
Further, 
\[\| \theta^n \| \leq \dt \max_{t \in [t^n,t^{n+1}]} \norm{d_t^2 u }_{L^2(\Omega)} \]
so that, using lemma \ref{lem:est_xipi} (projection error) and $\dt \le \kappa h$ and with $\taueps = 1+\epsinv$
\[\max_n  \abs{c^n} \leq \taueps \left[ C h \max_{t\in [0,T]} \norm{u(t,\cdot)}^2_{H^1(\Omega)} + (T\dt^2+C\dt^3) \max_{t \in [0,T]} \norm{d_t^2 u }^2_{L^2(\Omega)}\right].\]
This implies
\begin{align*}
    \norm{\xi_h^M}^2 + \sum_{j=0}^{M-1} \varepsilon\dt
    \seminormbeta{\xi_h^j}^2
    &\le T  (e^{\frac{\dt M}{T} }-1)\taueps
    \left[ C h
    \norm{u}^2_{C^0([0,T],H^1(\Omega))} + (T\dt^2+C\dt^3) \norm{u}^2_{C^2([0,T],L^2(\Omega))}\right]\\
    &\le T  (e^{\frac{\dt M}{T} }-1)\taueps
    \left[ Ch
    \norm{u}^2_{C^0([0,T],H^1(\Omega))} + T\dt^2 \norm{u}^2_{C^2([0,T],L^2(\Omega))}\right].
\end{align*}
In the last line, we dropped the higher order terms, assuming that $h$ and $\dt$ are sufficiently small.
Now we use the triangle inequality
\[
\norm{u^n-u_h^n} \le \norm{u^n-\Pi_h u^n} + \norm{\Pi_h u^n-u_h^n} = \norm{\xi_h^n} + \norm{\xi_{\pi}^n},
\]
the fact that $(\norm{\xi_h^n} + \norm{\xi_{\pi}^n})^2 \le 2\norm{\xi_h^n}^2 + 2 \norm{\xi_{\pi}^n}^2$,
and lemma \ref{lem:est_xipi} and $\dt M \le T$ to get the final result
\begin{align*}
    \norm{u^M-u_h^M}^2 &+ \sum_{j=0}^{M-1} \varepsilon\dt \seminormbeta{u^j - u_h^j}^2
    \leq 
    2\norm{\xi_h^M}^2 + 2 \norm{\xi_\pi^M}^2 + \sum_{j=0}^{M-1} 2\varepsilon \dt
    \left( \seminormbeta{\xi_h^j}^2
    + \seminormbeta{\xi_\pi^j}^2\right) \\
    &\le C(\varepsilon T+1) h \norm{u}^2_{C^0([0,T],H^1(\Omega))} \\
    &\qquad + T  (e^{\frac{\dt M}{T} }-1)\taueps
    \left[  Ch
    \norm{u}^2_{C^0([0,T],H^1(\Omega))} + T\dt^2 \norm{u}^2_{C^2([0,T],L^2(\Omega))}\right]\\
    &\le (C +  T e^{\frac{\dt M}{T} }C) \taueps h
    \norm{u}^2_{C^0([0,T],H^1(\Omega))} + \taueps T^2(e^{\frac{\dt M}{T} }-1)\dt^2 \norm{u}^2_{C^2([0,T],L^2(\Omega))}.
\end{align*}
\end{proof}

\section{Numerical results}\label{sec: numerical results}

In this section we present numerical results to examine the optimality of the theory established above.
We consider the ramp geometry depicted in figure \ref{fig: ramp geom a} for a variety of angles $\gamma$ ($\gamma=5^{\circ}, 15^{\circ},25^{\circ}, 35^{\circ},45^{\circ}$ degrees), together with the vector field
\[
\beta(x, y) = \frac{2 + \sin(\gamma) (x - x_0) - \cos(\gamma) y}{2 \sqrt{1 + \tan(\gamma)^2}} \begin{pmatrix}
    1\\
    \tan(\gamma)
\end{pmatrix}
\]
and initial data
\[
u_0(x, y) = \sin\left(\frac{\sqrt{2} \pi \big(\cos (\gamma) (x - x_0) + \sin(\gamma) y \big)}{1 - x_0}\right).
\]
Both depend on the choice of $\gamma$.
The point $(x_0, 0) = (0.2001, 0)$ marks the start of the ramp and the initial background mesh $\widehat{\mathcal{T}}_h$ covers the domain $(0,1)^2$.
The final simulation time is chosen as $T=0.5$. On the inflow boundary we prescribe as boundary data the analytic solution. The parameter $\tau$ in the definition of the capacity $\alpha_E$ in \eqref{eq:newdefalpha} is chosen as $\tau=1$ based on considerations of using $\tau \approx \frac{1}{\norm{\beta}_{\infty}}$ implied by lemma \ref{lem: inv est const}.

\begin{figure}[!htb]
    \centering
    \includegraphics[width=1.0\linewidth, height=0.4\linewidth]{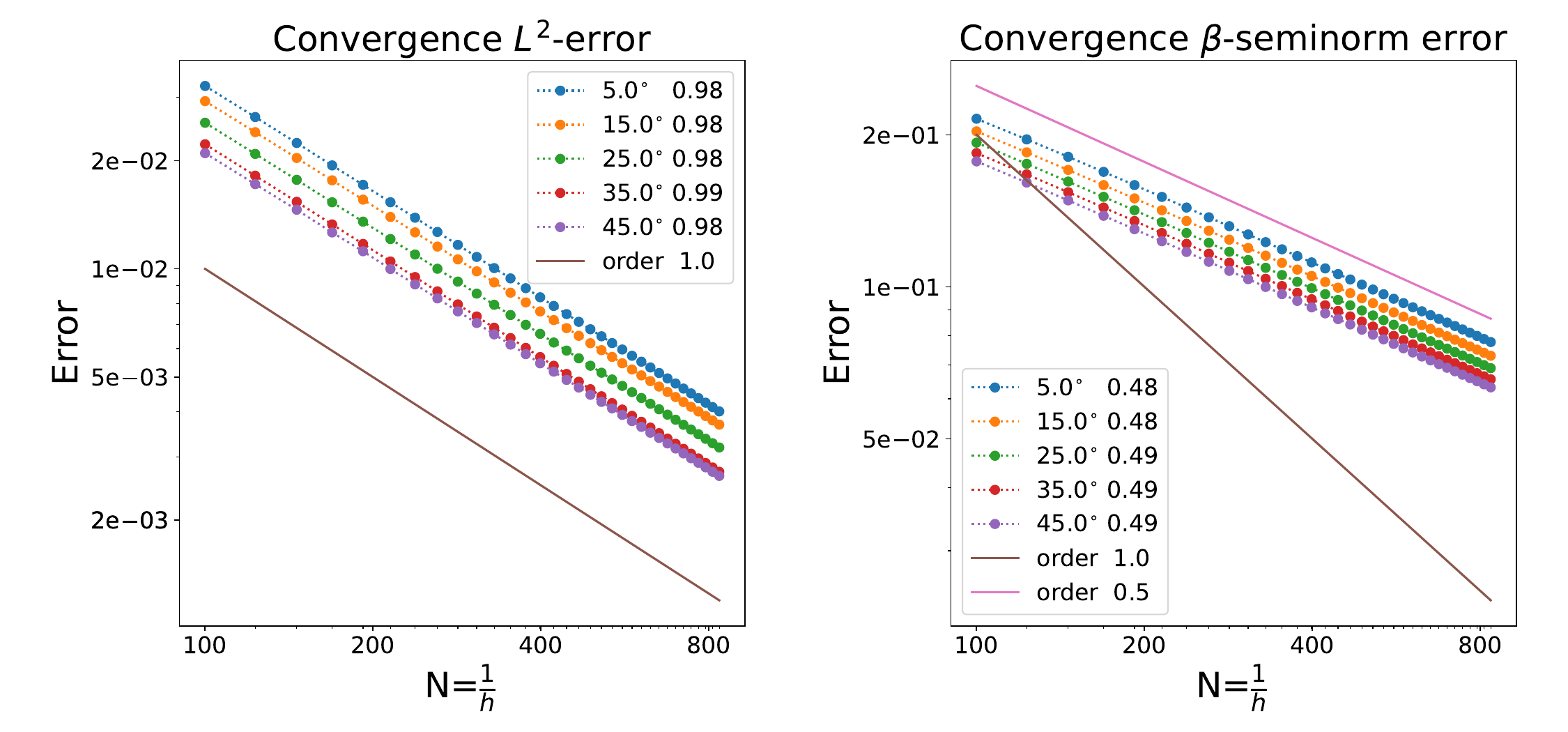}
    \caption{Convergence results for the time step choice $\dt = \frac{1}{5} \frac{h}{ \| \beta \|_{\infty} }$ given in \eqref{eq:theorem cfl}. \textit{Left:} Error at time $T$ measured in the $L^2$ (in space) norm. \textit{Right:} Error at time $T$ measured in the $\beta$-seminorm (in space).  The left plot shows convergence orders of $1$ for the error measured in the $L^2$ norm which is better than the order of $0.5$ predicted by theorem \ref{thm:main}. The right plot shows convergence orders of $0.5$ for the error measured in the $\beta$-seminorm at the final time $T$ which coincides with the predicted order of theorem \ref{thm:main} for the time-averaged $\beta$-seminorm.} 
    \label{fig:small-cfl-error}
\end{figure}


\begin{figure}[!htb]
    \centering
    \includegraphics[width=1.0\linewidth, height=0.4\linewidth]{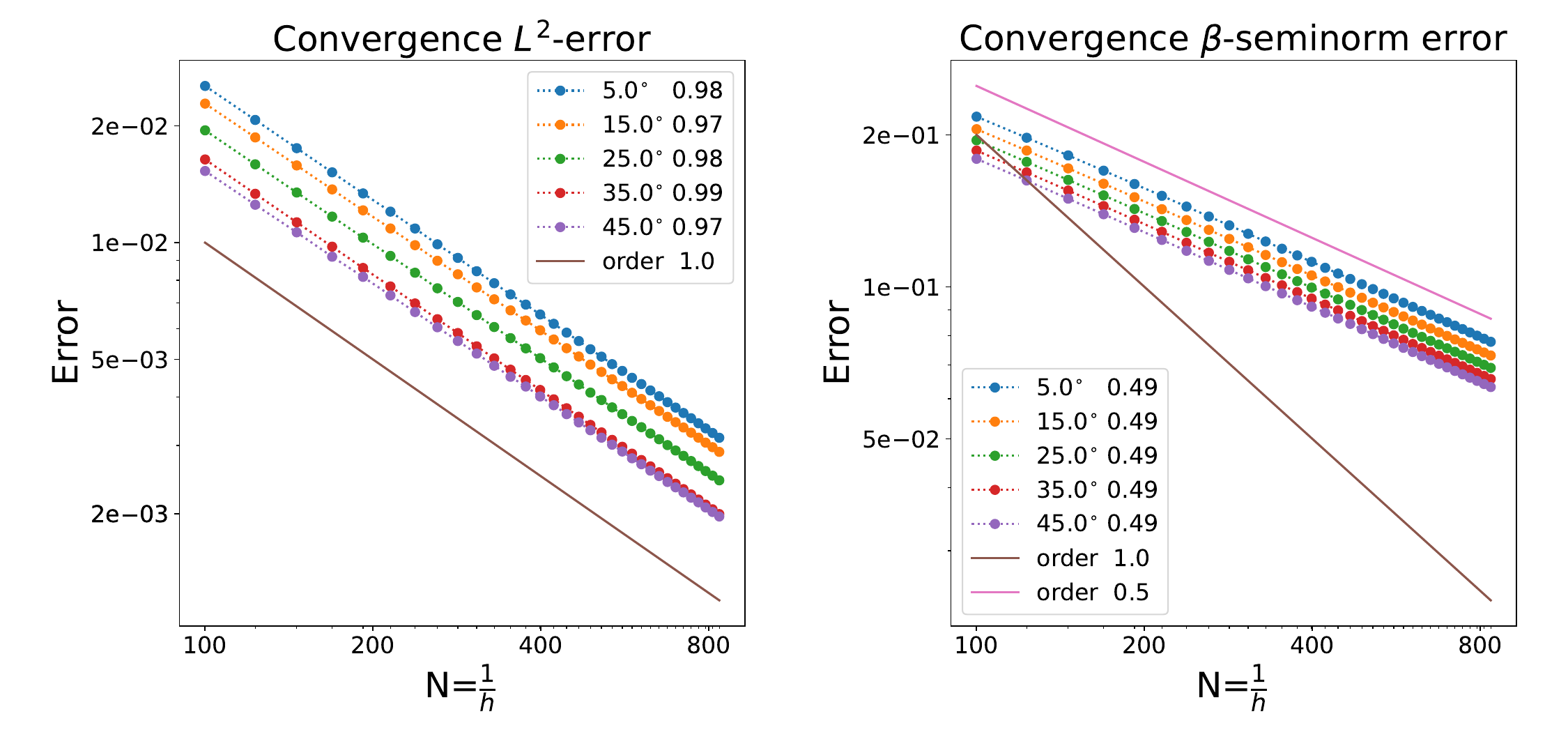}
    \caption{Convergence results for the time step choice $\dt = \frac{1}{2} \frac{h}{ \| \beta \|_{\infty} }$ given in \eqref{eq:usual cfl}. \textit{Left:} Error at time $T$ measured in the $L^2$ (in space) norm. \textit{Right:} Error at time $T$ measured in the $\beta$-seminorm (in space). Also outside of the provable regime we observe the same convergence order of 1 in the $L^2$ norm and of 0.5 in the $\beta$-seminorm, both at the final time $T$, as for the more restrictive time step choice.
    }
    \label{fig:usual-cfl-error}
\end{figure}

Choosing $\varepsilon$ in theorem \refwithname{thm:discrete stability} close to zero, our main result in theorem \ref{thm:main} requires a CFL condition of the form
\[
    \dt = c \frac{h} { \| \beta \|_{\infty}},\qquad\text{with $c < \frac{1}{4}$.}
\]
We will present numerical results for this time step constraint with $c = \frac{1}{5}$, i.e.
\begin{equation} \label{eq:theorem cfl}
    \dt = \frac{1}{5} \frac{h}{ \| \beta \|_{\infty} },
\end{equation}
for completeness. In practice, a bigger time step can be taken and we will further present numerical results for the time step constraint
\begin{equation} \label{eq:usual cfl}
    \dt = \frac{1}{2}\frac{h}{ \| \beta \|_{\infty}}.
\end{equation}


In figure \ref{fig:small-cfl-error} we present convergence results for the theorem-based time step choice given by \eqref{eq:theorem cfl}. We measure the error at the final time $T$ in both the $L^2(\Omega)$-norm and the $\seminormbeta{\cdot}$-seminorm in space. For all considered angles $\gamma$, we observe convergence of order $1$ for the $L^2$-norm in space and of order 0.5 for the $\seminormbeta{\cdot}$-seminorm in space. Our main result in theorem \ref{thm:main} implies convergence of order 0.5 for the norm used there. This indicates that our estimates are quasi-optimal for the norm of the error that is estimated but that the induced estimate for the error in the $\max_{0 \le n \le M} \norm{u(t^{n},\cdot)-u_h^n}_{L^2(\Omega)}$-norm is sub-optimal. It is unclear how to improve the estimate for this norm.

Generally, our presented error estimate is consistent with the result for quasi-uniform meshes in \cite{DiPietro_Ern}. That means, despite the presence of cut cells, we obtain the same convergence orders, measured in a slightly different norm as we modified the $\beta$-seminorm to account for the presence of cut cells.

Finally, in practical applications on quasi-uniform meshes larger time steps are used, without reducing the order of convergence. We investigated whether this also holds for the stabilized cut-cell formulation, using a time step size as given in
\eqref{eq:usual cfl}. Figure \ref{fig:usual-cfl-error} plots the errors and  the observed  convergence orders are almost identical to the convergence orders for the more restrictive CFL condition, showing again orders of 1 for the $L^2$-norm in space and of 0.5 for the $\beta$-seminorm. This means that also beyond the proven properties the stabilized method shows comparable properties as the original DG method.

\subsection*{Data availability}

Our numerical results are calculated with a C++ code that we based on the Dune \cite{dune-recent,dunepaperI:08,dunepaperII:08,dune-functions-1,dune-functions-2} libraries as the numerical core
and the TPMC \cite{tpmc} libraries to calculate the cut cell geometries.

The C++ source code to reproduce the presented numerical results is made available in \cite{zenodo-code}. 

\subsection*{Conflict of Interest}

The authors have no competing interests to declare that are relevant to the content of this
article.

\subsection*{Funding}
The research of JG, CE and GB was supported by the Deutsche Forschungsgemeinschaft (DFG, German Research Foundation) - SPP 2410 Hyperbolic Balance Laws in
Fluid Mechanics: Complexity, Scales, Randomness (CoScaRa), specifically JG within the project 525877563 (A posteriori error estimators for statistical solutions
of barotropic Navier-Stokes equations) and CE and GB within project 526031774 (EsCUT: Entropy-stable high-order CUT-cell discontinuous Galerkin methods).
SM, CE and GB gratefully acknowledges support by the Deutsche Forschungsgemeinschaft (DFG, German Research Foundation) - 439956613 (Hypercut).
CE and GB further acknowledge support by the Deutsche Forschungsgemeinschaft (DFG, German Research Foundation) under Germany's Excellence Strategy EXC 2044 –390685587, Mathematics Münster: Dynamics–Geometry–Structure.
JG also acknowledges support by the German Science Foundation (DFG) via
grant TRR 154 (Mathematical modelling, simulation and optimization using the example of gas
networks), sub-project C05 (Project 239904186).

\appendix

\bibliographystyle{alpha}
\bibliography{biblio}

\end{document}